\documentclass[12pt]{amsart}
\usepackage{graphicx}
\usepackage{amssymb}
\usepackage{setspace}
\usepackage{hyperref}
\usepackage{color}

\newtheorem{lemma}{Lemma}[section]
\newtheorem{proposition}[lemma]{Proposition}

\newtheorem{theorem}[lemma]{Theorem}

\newtheorem{corollary}[lemma]{Corollary}
\newtheorem{definition}[lemma]{Definition}

\newcommand{\R}{\mathbb{R}}
\newcommand{\C}{\mathbb{C}}
\newcommand{\h}{\mathbb{H}}
\newcommand{\N}{\mathbb{N}}
\newcommand{\Z}{\mathbb{Z}}

\newcommand{\RS}{\hat{\mathbb{C}}}
\newcommand{\PSL}{ {\rm PSL}({\rm 2},\mathbb{C})}
\newcommand{\psl}{ {\rm PSL}({\rm 2},\mathbb{R})}
\newcommand{\tli}{\tilde{l}_i}
\newcommand{\tl}{\tilde{l}}

\newcommand{\ol}{\overline}
\newcommand{\td}{\tilde}
\newcommand{\mc}{\mathcal}
\newcommand{\mb}{\mathbb}
\newcommand{\mr}{\mathring}

\newcommand{\sm}{\setminus}
\newcommand{\st}{\subset}

\newcommand{\pt}{\partial}

\newcommand{\cn}{\colon}

\newcommand{\ap}{\alpha}
\newcommand{\bt}{\beta}
\newcommand{\gm}{\gamma}
\newcommand{\ep}{\epsilon}

\newcommand{\oD}{\mathring{\mb{D}}^2}
\newcommand{\D}{\mb{D}}
\newcommand{\kp}{\kappa}

 \title[Decomposition of Projective Structures]{A Schottky Decomposition Theorem for Complex Projective Structures}
\author{Shinpei Baba}

\address{Mathematisches Institut, Universit\"{a}t Bonn, 53115 Bonn, Germany   }
\email{baba@math.uni-bonn.de}
\date{\today}

\begin{document}
\maketitle

\begin{abstract}
Let $S$ be a closed orientable surface of genus at least two, and let $C$ be an arbitrary (complex) projective structure on $S$.
We show that there is a decomposition of $S$ into pairs of pants and cylinders such that the restriction of $C$ to each component has an injective developing map and a discrete and faithful holonomy representation. 
This decomposition implies that every projective structure can be obtained by the construction of Gallo, Kapovich, and Marden.
Along the way, we show that there is an admissible loop on $(S, C)$, along which a grafting can be done. 
	
\end{abstract}

\section{Introduction}
Let $F$ be a connected orientable $C^1$-smooth surface possibly with boundary, and let $\td{F}$ denote the universal cover of $F$.
A {\it (complex) projective structure} $C$ on $F$ is a $(\RS, \PSL)$-structure, where $\RS = \C \, \cup \, \{\infty\}$ is the Riemann sphere.
In other words, it is a maximal atlas of $F$ modeled on $\RS$ with transition maps  in $\PSL$. 
The pair $(F, C)$ is called {\it a projective surface}.
As usual, we will often conflate $C$ and $(F, C)$.

There is an equivalent definition,  which we will mostly use in this paper:
A projective structure is a pair $(f, \rho)$, where $f: \tilde{F} \to \RS$ is a $C^1$-smooth locally injective map and $\rho: \pi_1(F) \to \PSL$ is a homomorphism, such that $f$ is $\rho$-equivariant, 
i.e. $f \circ \gamma = \rho (\gamma) \circ f$ for all $\gamma \in \pi_1(F)$.  
Then $f$ is called the {\it developing map} and $\rho$ the {\it holonomy (representation)} of the projective structure. 
On the interior of $\td{F}$, $~f$ is a local homeomorphism, and the restriction of $f$ to each boundary component of $\td{F}$ is a $C^1$-smooth curve.

A projective structure $C = (f, \rho)$ is defined up to an isotopy of $F$ and the action of an element of $\PSL$, i.e. the post-composition of $f$ with $\gamma \in \PSL$  and the conjugation of $\rho$ by $\gamma$.
(See \cite[3.4]{Th2}, \cite[7.1]{K2}.)
The $C^1$-smoothness is required to define a natural topology on the space of all projective structures on $F$ in the case that $F$ is not closed (see \cite[I.1.5]{CEG}). 
In this paper, we always assume that surfaces are connected and $C^1$-smooth possibly with boundary (although we do not consider the deformation space of projective structures). 

\begin{definition}
A projective structure $C = (f, \rho)$ on $F$ is called {\bf admissible} 
if it suffices the following conditions: \\
\noindent (i) $f$ is an embedding and \\
\noindent (ii)  $\rho$ is an isomorphism onto a quasifuchsian group (if $F$ is closed) or  onto a Schottky group (if $F$ is not closed). \\  
Analogously, a simple loop $l$ on the projective surface $(F, C)$ is {\bf admissible} if $\td{l}$ embeds in $\RS$ by $f$ and $\rho(\gamma_l)$ is loxodromic, where $\td{l}$ is a lift of $l$ to $\td{F}$ and $\gamma_l$ is the homotopy class of $l$. 
\end{definition}
\noindent \emph{Remark:}  Condition (ii) is equivalent to saying that $\rho$ is an isomorphism onto a convex-cocompact subgroup of $\PSL$.

A hyperbolic structure is a basic example of a projective structure, since  $\h^2 \st \RS$ and ${\rm Isom}^+(\h^2) \cong {\rm PSL}(2, \R) \st \PSL$ in a compatible way.  
Every hyperbolic structure on a closed orientable surface is an admissible projective structure.
However, in general, developing maps are not necessarily injective and holonomy representations are not necessarily discrete or faithful (c.f. (iii) and (iv) following Corollary  \ref{Pants} in this section). 
Throughout this paper, let $S$ denote a closed orientable surface of genus at least $2$ and let $\td{S}$ denote the universal cover of $S$.
(The orientability of $S$ is not essential for the mains theorems of this paper, if we consider two-dimensional Mobius structures instead of projective structures.)
The following theorem yields a decomposition of an arbitrary projective surface into admissible projective subsurfaces:\\

\noindent \textbf{ Theorem \ref{id}}
Let $C$ be a projective structure on $S$. 
Then there exists a decomposition of $S$ into  cylinders and compact connected surfaces of negative Euler characteristic, such that the restriction of $C$ to each cylinder is an integral flat structure and the restriction to each surface of negative Euler characteristic is an admissible projective structure.  \\

By a \emph{decomposition}, we mean that the subsurfaces in this theorem are the connected components of $S$ minus some \emph{multiloop}, disjoint union of essential simple loops.  

An {\it integral flat structure} is a basic projective (actually affine) structure on a cylinder, closely related to an operation called {\it grafting} (\S \ref{FlatCylinders}).
If there is  an admissible loop on a projective surface, we can define a grafting along this loop  (see  \cite{K2}, \cite{Go2}, \cite{Bromberg} for example). 
This operation gives another projective structure on the same surface, preserving the orientation and the holonomy representation.  
If the admissible loop is circular (i.e. it corresponds to a simple circular arc on $\RS$ via the developing map),  
then the integral flat structure is exactly the projective structure that the grafting operation inserts to the projective structure along the loop. 
An integral flat structure on an annulus can be easily decomposed into integral flat structures that are admissible.

Theorem \ref{id} immediately implies:\\

\noindent \textbf{Corollary \ref{Pants}}
Let $C$ be an arbitrary projective structure on $S$.
Then there exists a decomposition of $S$ into pairs of pants and cylinders such that the restriction of $C$ to each cylinder is an integral flat structure and the restriction to each pair of pants is an admissible structure.  \\
 
This corollary gives the affirmative answer to a question raised by Gallo, Kapovich, and Marden (\cite[12.1]{GKM}). 
The authors of \cite{GKM} gave necessary and sufficient conditions for a representation $\rho: \pi_1(S) \to \PSL$ to be the holonomy representation of some projective structure on $S$.
The conditions are:  (iii) $Im(\rho)$ is a non-elementary subgroup of $\PSL$ and (iv) $\rho$ lifts to a representation from $\pi_1(S)$ to ${\rm SL}({\rm 2}, \C)$. 
In order to prove the sufficiency of these conditions, given an arbitrary representation $\rho$ satisfying (iii) and (iv), they constructed a projective structure on $S$ with holonomy $\rho$ in the following way:
First, decompose $S$ into pairs of pants, $\{P_i\}$, such that $\rho|_{\pi_1(P_i)}$ is an isomorphism onto a rank-two Schottky group for each $i$.
Second, construct an admissible projective structure on each $P_i$ with the holonomy representation $\rho|_{\pi_1(P_i)}$.
Last, glue these structures on the pairs of pants together by inserting projective structures on cylinders between the corresponding  boundaries of $P_i$'s, and obtain a desired projective structure.
They asked whether every projective structure on $S$ arises from such a Schottky pants decomposition. 
More specifically, they asked if every projective structure contains an admissible loop, which is answered by: \\

\noindent \textbf{Theorem  \ref{gc}}
\ For every projective structure $C$ on $S$, there exists an admissible loop on $(S, C)$.\\

\emph{A remark on Theorem \ref{gc}:}
From our argument, it immediately follows that, on every projective surface $(S, C)$, there are infinitely many homotopy classes of admissible loops, unless $(S,C)$ satisfies the following condition (which almost always fails):  
\begin{itemize}
\item Setting  $C = (\tau, L)$ to be the expression of $C$ in Thurston's coordinates (\S \ref{ThP}), $L$  decomposes $S$ into a disjoint union of pairs of pants.
\end{itemize}
The author believes that, on every $(S,C)$ (even with this condition), each simple loop sufficiently close to $L$, in $\mathcal{PML}(S) \cong \mb{S}^{6g-7}$, is homotopic to an admissible loop.\\ %

Note that Theorem \ref{gc} is weaker than Theorem \ref{id}, since the boundary components of each subsurface in Theorem \ref{id} are, in particular, admissible loops. 
However, the proof of Theorem \ref{gc} contains the basic ideas for the proof of Theorem  \ref{id}.
Furthermore, Theorem \ref{gc} addresses the following question about grafting: \\

\noindent \textbf{Question 1} (\cite[\S 12]{GKM}): Assume that two projective structures on $S$ have the same orientation and holonomy representation. Can one projective structure be transformed to the other by a sequence of graftings and inverse-graftings?\\

The grafting and inverse-grafting operations generate an equivalence relation among the projective structures with a given holonomy representation. 
Question 1 asks if there are exactly two equivalence classes represented by the orientations of the projective structures.   
Theorem  \ref{gc}  implies that every equivalence class consists of infinitely many projective structures. 
Ultimately, Question 1 aims to characterize the collection of projective structures with the given holonomy representation; this characterization problem goes back, at least, to Hubbard's  paper published in 1981 (\cite{Hu}; see also \cite{K1}, \cite{GKM}). 
In the special case that the holonomy representation is an isomorphism onto a  quasifuchsian group, the characterization is given by Goldman, using grafting, and the answer to Question 1 is affirmative (\cite[Theorem C]{Go1}).

The {\it holonomy map} $Hol: P(S) \to V(S)$ is a projection given by $C = (f, \rho) \mapsto \rho$, where $P(S)$ is the space of all projective structures on $S$ and $V(S)$ is the representation variety of homomorphisms from $\pi_1(S)$ to $\PSL$.
This map is {\it not} a covering map onto its image (\cite{Hej}), which makes problems in this area difficult. 

One may ask  the above questions in the case of other $(G,X)$-structures as well (c.f.  \cite[1.10]{Go1}). 
In particular,  S. Choi gave a canonical decomposition of real projective structures, i.e. $({\rm PGL}(3, \R), \R {\rm P}^2)$-structures, analogous to the one given by Theorem \ref{id} (see \cite{C1}).  \\

\noindent \textbf{ An outline of the proofs:}
For a given projective surface $(S, C)$, there is a corresponding pair $(\tau, L)$ of a marked hyperbolic structure $\tau$ on $S$ and a measured geodesic lamination $L = (\lambda, \mu)$ on $(S, \tau)$ (see Thurston's coordinates in \S \ref{ThP}).
A periodic leaf of $\lambda$ corresponds to a continuous family of admissible loops on $(S, C)$.
For each irrational minimal sublamination of $\lambda$, consider a standard sequence $(l_i)$ of simple geodesic loops on $(S, \tau)$ that approximates $\lambda$ (\S \ref{Appro}). 
Then we will show that $l_i$ is admissible for all sufficiently large $i$.
These admissible loops yield Theorem \ref{gc}.
Taking a disjoint union of such admissible loops, we construct a multi-loop on $(S, C)$ that is a good approximation of the entire lamination $\lambda$. 
We will show that the multi-loop on $(S, C)$ achieves the desired admissible decomposition described in Theorem \ref{id}.

Every measured lamination on $\h^2$ induces a continuous map from $\h^2$ to $\h^3$, called a {\it bending map} (\S \ref{bending}).
Via the bending map, the measured lamination corresponds to a projective structure on $\oD$ (\S \ref{ThP}).
Our proofs are based on the fact that injective quasiisometric bending maps correspond to admissible projective structures on $\oD$ (\S \ref{IBM}).  
In order to show that $l_i$ is admissible, we take the {\it total lift} $\td{L}$ of $L$ to $\h^2$ (\S \ref{DT}) and a lift $\td{l}_i$ of $l_i$ to $\h^2$.
Consider the sublamination $I(\td{L}, \tli)$ of $\td{L}$ that consists of the leaves of $\td{L}$ intersecting $\td{l}_i$, so that the structure on $\tli$ embeds into the projective structure on $\oD$ corresponding to $I(\td{L}, \tli)$. 
In other words,
$\td{L}$ and $I(\td{L}, \tli)$ coincide in a sufficiently small neighborhood of $\tli$, and therefore, $I(\td{L}, \tli)$ is sufficient to capture the structure on $\tli$.  
The transversal measure of  $I(\td{L}, \tli)$ is uniformly small (Proposition \ref{propA}).
Then, accordingly, the bending map  induced by $I(\td{L}, \tli)$ bends $\h^2$ inside $\h^3$ to a uniformly small degree,  and  therefore it is an injective quasiisometric embedding (\S \ref{IBM}).  
Therefore $I(\td{L}, \tli)$ corresponds to an admissible structure on $\oD$. 
Since the structure on $\td{l}_i$ is embedded in the admissible structure on $\oD$,  $~l_i$ is also admissible. 
Theorem \ref{id} will be proven based on the same idea.

\textbf{Acknowledgments.} The author would like to thank my advisor, Misha Kapovich, for introducing me to complex projective structures and supporting me with his comments and encouragements.  
He also thanks the referee for valuable comments.
The author was partially supported by NSF grants.

\section{Notation}

\noindent $\ol{xy}$: the geodesic segment connecting $x$ and $y$, where $x$ and $y$ are points in a metric space.\\
$cl(X)$: the closure of $X$, where $X$ is a subset of a topological space.\\
$\mr{X}$: the interior of $X$.\\
$Conv(Y)$: the convex hull of $Y$, where $Y$ is a subset of  $\ol{\h^3}$.\\
$D_r(x)$: the closed disk of radius $r$ centered at $x$ in a hyperbolic space.\\
\section{Preliminaries}\label{P}

\subsection{Measured Laminations}\label{ml}
(For details, see \cite{PH}, \cite{CEG}, \cite{CB}, \cite{K2},  \cite{Th2}.)
Let $F$ be a Riemannian surface with a constant curvature, possibly with geodesic boundary. 
A {\it geodesic lamination} on $F$ is a collection of disjoint simple complete geodesics on $F$ whose union is a closed subset of $F$.
Each geodesic of a geodesic lamination is called a {\it leaf}. 
For a measured lamination $\lambda$  on $F$, let $| \lambda | \>   \st F $ denote the union leaves of $\lambda$.
If $|\lambda| = F$, then $\lambda$ is called a {\it foliation}.
A {\it measured lamination} is a pair  $L = (\lambda, \mu)$, where $\lambda$ is a geodesic lamination, and $\mu$ is a transversal measure of $\lambda$.
Let a leaf of $L$ refer to a leaf of $\lambda$ and $| L |$ refer to $| \lambda |$.
We always assume that $\lambda$ is the support of $\mu$, i.e.  if $l$ is a leaf of $L$,  then $\mu(s) > 0$ for every geodesic segment $s$ that transversally intersects $l$.  
The {\it weight} of a leaf $l$ of $L$ is $ \inf\{ \mu(s) \}$, where $s$ varies over all geodesic segments that transversally intersect $l$, i.e. the atomic transversal measure of $l$.
We denote the weight of $l$ by $w(l)$.
If $l$ is an isolated leaf of $L$, then $w(l) = \mu(s)$ for every geodesic segment $s$ on $F$  that transversally intersects $|L|$ exactly once at a point on $l$.

By convention, if a geodesic segment $s$ is contained in a leaf of $\lambda$, then $\mu(s) = 0$.
In addition, letting $x, y$ be the end points of $s$, if  $x$ or $y$ is contained in a leaf with positive weight, then the weight does not contribute to the value of $\mu(s)$, so that $\mu(s) = \sup \mu(s')$ where $s'$ varies over all geodesic segments strictly contained in $s$, i.e.  $s' \st s \sm \{x, y \}$.

Recall that $S$ is a closed orientable hyperbolic surface.
A measured lamination on $S$ is {\it minimal} if it does not contain any proper sublamination. 
Every measured lamination $L$ on $S$ uniquely decomposes into a finite number of disjoint minimal laminations of the following two types:
a periodic leaf with positive weight ({\it periodic minimal lamination}), and  a measured lamination consisting of uncountably many bi-infinite geodesics ({\it irrational minimal lamination}).

Let $M = (\nu, \omega)$ be an irrational minimal sublamination of $L$.
Then, since each minimal sublamination of $L$ forms a closed subset of $S$,  $~|\nu|$ is also an open subset of $|\lambda|$.
Here are other properties of $M$:
Each leaf of $\nu$ is a dense subset of $|\nu|$.
If $s$ is a geodesic segment on $S$ transversally intersecting $\nu$, then $s \cap |\nu|$ is a Cantor set  i.e., a closed, perfect subset with empty interior (\cite[Corollary 1.7.6, 1.7.7]{PH}). 
Therefore, no leaf of $M$ has a positive weight. 


\subsection{Flat Cylinders}\label{FlatCylinders}

For every $\theta \in (0, 2\pi)$ and distinct $z_1, z_2 \in \RS$, let  $R_\theta$  be the open region in $\RS$ bounded by two simple circular arcs connecting $z_1$ and $z_2$ such that the inner angles at the vertices $z_1$ and $z_2$ are equal to $\theta$.  
Since $R_\theta$ is embedded in $\RS$, it is equipped with a canonical projective structure (whose developing map is the identity map).
We call the structure a {\it crescent} of angle $\theta$; its projective structure only depends on the choice of $\theta$.

Let $\alpha$ be a hyperbolic element of $\PSL$ that fixes $z_1$ and $z_2$.
Then $\langle \alpha \rangle$, the subgroup of $\PSL$ generated by $\alpha$, is an infinite cyclic group acting on $R_\theta$ freely and property discontinuously.
By quotienting $R_\theta$ by $\langle \alpha \rangle$, we obtain a projective structure on a cylinder.
We call the structure a {\it flat structure} of height $\theta$; it only depends on $\theta$ and the translation length of $\alpha$.
A flat structure of height $\theta$ on  a cylinder forms a projective surface, which we call a {\it flat cylinder} of height $\theta$.

We shall define a crescent and a flat cylinder for arbitrary $\theta > 0$, generalizing those for  $0 < \theta < 2\pi$.  
For arbitrary $\theta > 0$, let $R_\theta = (0,\theta) \times (0, \infty) \st \R^2$.
Define $f_\theta: (0, \theta) \times (0, \infty) \to \C$  by $f_\theta(x, y) = y \cos x + \sqrt{-1} ~y \sin x$ (i.e. the polar coordinates). 
Then $f_\theta$ defines a projective structure on $R_\theta \cong \oD$, which is a crescent of angle $\theta$.
For $a > 0$, let $T_a$ be an automorphism of $\R^2$ defined by $T_a(x,y) = (x, a y)$.
Define a homomorphism $\rho:  \langle T_a \rangle \cong \Z \to \PSL$ by $\rho(T_a)(z) = a z$ for all $z$.
Then $f_\theta$ is $\rho$-equivariant. 
Quotienting $R_\theta$ by the action of $\langle T_a \rangle$, we obtain a flat cylinder of height $\theta$.
A flat cylinder of height $\theta$ is {\it integral} if $\theta$ is a multiple of $2\pi$.  
The multiplier is called the {\it degree} of the integral flat cylinder, so that, for all $z \in \C \sm \R_{\geq 0}$, the cardinality of $f_\theta^{-1}(z)$ is equal to the degree.  
Clearly, an integral flat cylinder of height $n$ can be decomposed into $n$ flat cylinders of height one, which are admissible.  

Let $C = (f_\theta, \rho_{id})$ be the crescent of angle $\theta > 0$ given in the form above, where $\rho_{id}: \pi_1(R_\theta) \to \PSL$ is the trivial representation.
For each $x \in (0, \theta)$,  $f_\theta$ takes $x \times (0,\infty)$ to a straight line on $\C$ connecting $0$ and $\infty$.
The collection of these lines,  $\{~ \{x\} \times (0, \infty)~ |~ x \in (0, \theta) \}$,  forms a foliation $\lambda_C$ on $R_\theta$, which we call the {\it canonical foliation} on $C$.
We also can define the canonical transversal measure $\mu_C$ of $\lambda_C$ by $\mu_C(\ol{P_1P_2}) = |x_1 - x_2|$ for all $P_1 = (x_1, y_1), P_2 = (x_2, y_2) \in R_\theta$. 
Note that $\{x_1\} \times (0, \infty)$ and $\{x_2\} \times (0, \infty)$ bound a crescent of height $\mu_C(\ol{P_1P_2})$ contained in $C$.
Call $(\lambda_C, \mu_C)$ the {\it canonical measured foliation} on $C$. 

For each $y \in (0, \infty)$,  $~(0, \theta) \times \{y\}$ is orthogonal to (each leaf of) $\lambda_C$ in terms of the angles obtained by pulling back the conformal structure on $\RS$ via $f_\theta$.
Besides, $f_\theta$ takes $(0, \theta) \times \{y\}$ to a (not necessarily simple) circular arc  on $\RS$.
The collection of these orthogonal lines  $\{ (0, \theta) \times \{y\} ~|~ y \in (0, \infty) \}$  forms a foliation on $R_\theta$, which is dual to $\lambda_C$.
By identifying the points on each leaf of the dual foliation, $R_\theta$ projects to a line.   
Since $(\lambda_C, \mu_C)$ and its dual foliation are invariant under the action of $\langle T_a \rangle$, we obtain the {\it canonical foliation} and its dual foliation on the flat cylinder $C/ \langle T_\alpha \rangle$. 
Accordingly, the flat cylinder projects to a circle by identifying the points on each leaf of the dual foliation.

Since  $f_\theta: (0, \theta) \times (0, \infty) \to \C$ continuously extends to $\{0, \theta\} \times (0, \infty)$,
we can compactify a flat cylinder of height $\theta$ to a projective structure on a compact cylinder with boundary.
By abusing the notation, we call this compactified flat cylinder, also, a flat cylinder of angle $\theta$.
Accordingly, the universal cover of the compactified flat cylinder of angle $\theta$ also is called a crescent of angle $\theta$.

\begin{figure}[htbp]
\includegraphics[width = 5in]{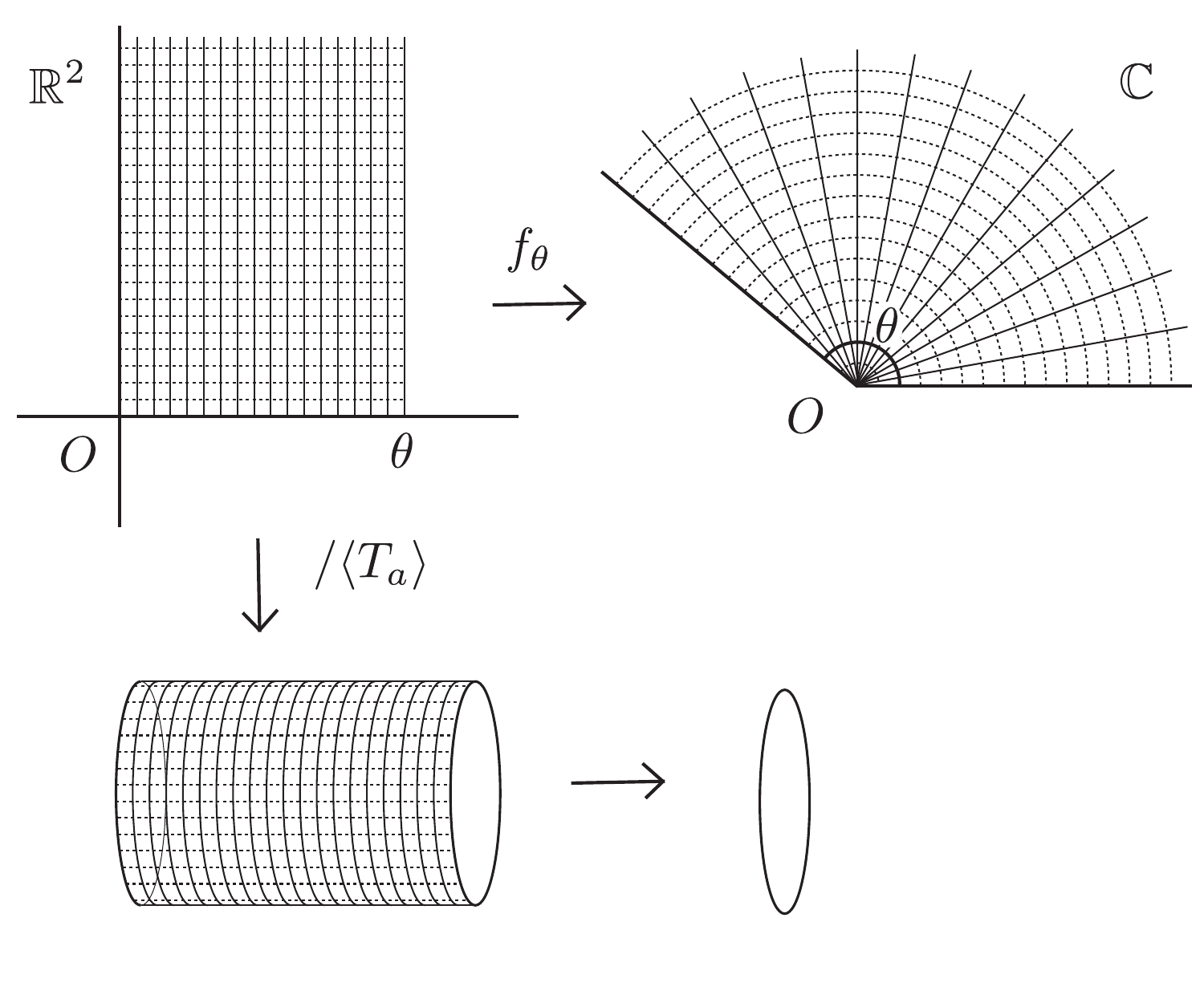}
\caption{}
\label{PicFlatCylinder}
\end{figure}

\subsection{Dual Tree}\label{DT} (For more details, see \cite{Morgan}, \cite{K2}.)
Let $S$ be a closed hyperbolic surface and $L = (\lambda, \mu)$ be a measured lamination on $S$ without periodic leaves.
Then no leaf of $L$ has a positive weight. 
Let $p: \h^2 \to S$ be the covering map.
The {\it total lift} of $L$ is a measured lamination $\td{L} = ( \td{\lambda}, \td{\mu})$ on $\h^2$, where $\td{\lambda}$ consists of all the lifts of the leaves of $\lambda$,  and $\td{\mu}$ is the pull back of $\mu$, so that $(\h^2, \td{L})$ is locally isomorphic to $(S, L)$ via $p$. 
Then $\td{L}$ is a $\pi_1(S)$-invariant measured lamination on $\h^2$. 

There is  a unique $\R$-tree dual to $\td{L}$ constructed in the following way (if $L$ contains periodic leaves, the construction is more complicated).
The transversal measure $\td{\mu}$ defines a pseudo-metric $d_{\td{\mu}}$ on $\h^2$,  by $d_{\td{\mu}}(x,y) = \td{\mu}( \ol{xy} )$ for all $x,y \in \h^2$. 
Since no leaf of $\td{L}$ has a positive weight, $d_{\td{\mu}}: \h^2 \times \h^2 \to \R_{\geq 0}$ is continuous. 
Define an equivalence relation on $\h^2$ by $x \sim y$ if and only if $\td{\mu}(\ol{xy}) = 0$.
There are only two types of equivalence classes:
the closure of a complementary region of $|\td{\lambda}|$,
and a leaf of $\td{\lambda}$ that is {\it not} a boundary geodesic of such a complementary region. 
Let $T$ be the quotient of $\h^2$ by the equivalence relation, and let $\mc{P}: \h^2 \to T$ be the quotient map.
In particular, if a geodesic in $\h^2$ is a leaf of $\td{L}$ or it is contained in the closure of a complementary region of $|\td{\lambda}|$, then $\mc{P}$ takes this geodesic to a point in $T$.
It turns out that $T$ is an $\R$-tree equipped with a canonical metric $d_T$ induced by $\td{\mu}$ (see \cite[\S 2.3]{O}):
For every $x, y \in T$, $~d_T (x, y) = \td{\mu} (\ol{x'y'})$, where $x', y' \in \h^2$ are such that $\mc{P}(x') = x$ and $\mc{P}(y') = y$.  
If a geodesic in $\h^2$ transversally intersects $\td{L}$, then $\mc{P}$ takes this geodesic to a geodesic in $T$. 
Since the action of $\pi_1(S)$ on $\h^2$ preserves $\td{L}$, it isometrically acts on $T$.
Let $l$ be a closed geodesic on $S$ transversally intersecting $L$, $~\td{l}$ be a lift of $l$ to $\h^2$, and $\gamma_l$ be the homotopy class of $l$ in $\pi_1(S)$. 
Then $\mc{P}(\td{l})$ is a geodesic in $T$, and the action of $\gamma_l$ on $T$ isometrically translates along $\mc{P}(\td{l})$ by the distance $\mu(l)$.
 
\subsection{Flow Boxes}

Let  $(a,b)$ and $(c,d)$ be open intervals in $\R$, and let $Y$ be a closed subset of $(c,d)$.
Consider a geodesic lamination $\lambda$ on $(a, b) \times (c,d) \st \R^2$ that consists of the leaves $(a, b) \times \{y\}$ for all $y \in Y$ (i.e. we have $Y$-worth of horizontal leaves). 
Let $\mu$ be a transversal measure for $\lambda$, and let $L$ be the measured lamination $(\lambda, \mu)$.  
The pair $((a,b) \times (c,d), L)$ is called a {\it(Euclidean) flow box}.
Let $s$ be a vertical geodesic in $(a,b) \times (c,d)$, i.e. $s = \{x\} \times (c,d)$ for some $x \in (a,b)$.
The {\it height} of the flow box is $\mu(s)$, which does not depend on the choice of $x$.  
If $L$ is a measured lamination on a hyperbolic quadrilateral $Q$, and $(Q, L)$ is isomorphic to a flow box, then we also call $(Q, L)$ a (hyperbolic) flow box. 

For every point $x$ on a leaf of $L$, there is a neighborhood of $x$ that is isomorphic to the a flow box (\cite[\S 11.6]{K2}, \cite[\S 1.6]{PH}); then we immediately obtain: 

\begin{lemma}\label{box-nhd}
Let $L = (\lambda, \mu)$ be a measured geodesic lamination on $\h^2$ without leaves of positive weight. 
Let $s$ be a geodesic segment contained in a leaf of $L$.
For every $\epsilon > 0$, there exists a neighborhood of $s$ isomorphic to a flow box of height less than $\epsilon$.
\end{lemma}
\noindent {\it Remark:} Under the projection map from $(\h^2, L)$ to its dual tree, the flow box neighborhood of height less than $\epsilon$ projects to a geodesic segment of length less than $\epsilon$.

\subsection{Approximating an Irrational Lamination}\label{Appro}
\begin{figure}[htbp]
\includegraphics[width= 5in]{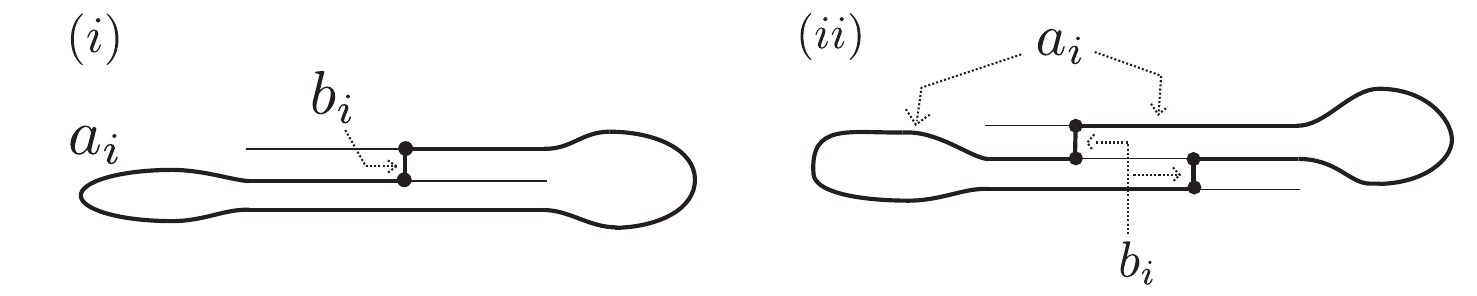}
\caption{}
\label{a_i&b_i}
\end{figure}

(For details, see \cite[I.4.2.15]{CEG}.)
Let $S$ be a closed (orientable) hyperbolic surface, and $L = (\lambda, \mu)$ be an irrational minimal measured lamination on $S$. 
Then there is a sequence $(l_i)$ of simple closed geodesics, limiting to $|L|$ in the Chabauty topology, where each $l_i$ is homotopic to a piecewise geodesic loop $c_i$ of one of the following two types: 
\begin{itemize}
 \item $c_i$ is a union of a (long) geodesic segment $a_i$ contained in a leaf of $L$ and a (short) geodesic segment $b_i$ transversal to $|L|$ (see Figure \ref{a_i&b_i} (i)).
\item  $c_i$ is a union of two (long) geodesic segments contained in a leaf of  $L$ and two (short) geodesic segments transversal to $|L|$; we let $a_i$ denote the union of these long geodesic segments, and $b_i$ denote the union of these short geodesic segments. (See Figure \ref{a_i&b_i} (ii).) 
\end{itemize}
In both cases, we have $\lim_{i \to \infty} length(a_i) \to \infty$ and $\lim_{i \to \infty} length(b_i) \to 0$.
Then, since $S$ is closed, $\lim_{i \to \infty} \mu(b_i) = 0$.

Let $m$ and $n$ be geodesics on a complete hyperbolic surface.
Assume that $m \cap n \neq \emptyset$, and pick  $p \in m \cap n$.
We can ``rotate'' $m$ to $n$ about $p$ by a unique angle in $(-\pi/2, \pi/2]$.
More precisely, we do the following:
Let $\td{p}$ be a lift of $p$ to $\h^2$,  and let $\td{m}$ and $\td{n}$ be the lifts of $m$ and $n$, respectively, to $\h^2$ that intersect at $\td{p}$.
Then, we can indeed rotate $\td{m}$ to $\td{n}$ about $\td{p}$ by a unique angle in $(-\pi/2, \pi/2]$.
Let $\angle_p(m,n)$ denote this angle,  and call it the {\it angle} between $m$ and $n$ at $p$. 

Let $\nu$ be a geodesic lamination on the surface that intersects $n$. 
\begin{definition}
The {\bf angle} between $\nu$ and $n$ is
$$ \angle(\nu, n) = \displaystyle \sup \{\, |\angle_p ( m, n)| ~\, | ~   m \in \nu ~{\rm and}~ p \in m \cap n~\} \in [0, \pi/2].$$
\end{definition}


\begin{lemma}\label{ml1}
$\displaystyle \lim_{i \to \infty}\angle(\lambda, l_i) = 0$ and $\displaystyle \lim_{i \to \infty}\mu(l_i) = 0$.
\end{lemma}

\begin{proof}  
Choose $x_i \in | \lambda | \cap l_i$, and let $m_i$ be the leaf of $\lambda$ that intersects $l_i$ at $x_i$.  
Let $\theta_i = \angle_{x_i }(m_i , l_i)$. 
Note that $|\lambda|$ is a compact subset of $S$.
Therefore, by taking a subsequence if necessary, $(x_i, \theta_i)$ converges to $ (x, \theta) \in |\lambda| \times [-\pi/2, \pi/2]$. 
Assume that $\displaystyle \lim_{i \to \infty}\angle(\lambda, l_i) \neq 0$. 
Then there exists a sequence $(x_i, \theta_i)$ converging to $(x, \theta)$ with $\theta \neq 0$.
The sequence $(l_i)$ converges to a geodesic that intersects $|\lambda|$ at $x$ with angle   $\theta$. 
This contradicts the convergence of $(l_i)$ to $|\lambda|$.

Next, we prove that $\lim_{i \to \infty} \mu(l_i) = 0$.
Let $c_i = a_i \cup b_i$, which is a simple loop on $S$.
Observe that $\lim_{i \to\infty}\mu(c_i) = 0$, since $\mu(a_i) \equiv 0$ and $\lim_{i \to \infty} \mu(b_i) = 0$.  
Therefore, it suffices to show that $\mu(l_i) \leq \mu(c_i)$.
The basic idea is that a geodesic loop realizes the minimal transversal measure among all the loops in the same homotopy class.
Let $\gamma_i \in \pi_1(S)$ be the homotopy class of $l_i$.
Regarding $l_i$ as a bi-infinite geodesic, let $\td{l}_i$ be the lift of $l_i$, such that $\gamma_i$ translates $\h^2$ along the geodesic $\td{l}_i$ by $length(l_i)$.
Regarding $l_i$ as a simple closed path, let $\bar{l}_i$ be a lift of $l_i$, such that that $\bar{l}_i$ is contained in $\td{l}_i$.
Then $\gamma_i$ identifies the ends of $\bar{l}_i$.
Recall that $\mc{P}$ is the projection from $\h^2$ to the $\R$-tree $T$ dual to $\td{L}$ (\S \ref{DT}). 
Then $~\mc{P}(l_i)$ is a geodesic in $T$, and the action of $\gamma_i$ on $T$ isometrically translates along $\mc{P}(l_i)$ by $\mu(l_i) = length(\mc{P}(\bar{l_i}))$.
Similarly, regarding $c_i$ as a simple closed path, let $\bar{c}_i$ be a lift of $c_i$ to $\h^2$, such that the ends of  $\bar{c}_i$ are identified by $\gamma_i$. 
Then $\mc{P}(\bar{c_i})$ is a piecewise geodesic path in $T$, and $\gamma_i$ identifies the ends of $\mc{P}(\bar{c_i})$.
We also have $\mu(c_i) = length(\mc{P}(\bar{c_i}))$.
The translation length of $\gamma_i$ is equal to or less than the distance between the ends of  $\mc{P}(\bar{c_i})$ since the translation length of $\gamma_i$ is $\inf \{ dist_T(x, \gamma_i(x)) \>|\>  x \in T \}$.
Since $\gamma_i$ translates along $P(\tli)$, the translation length is $length(\mc{P}(\bar{l}_i))$.
Hence, $\mu(l_i) \leq \mu(c_i)$.
\end{proof}

\subsection{Bending Maps} (For details, see \cite[II.1]{EM}, \cite{KT},  \cite{KP}.)\label{bending}
Let $L = (\lambda, \mu)$ be a measured lamination  on $\h^2$.
Then $L$ induces a {\it bending map} $\beta_L = \beta: \h^2 \to \h^3$ by ``bending $\h^2$ inside $\h^3$ along $\lambda$ by angle $\mu$".
The bending map $\beta_L$ is continuous and unique up to the post-composition with an element of $\PSL$.
In addition, $\beta_L$ is isometric on each leaf of $\lambda$ and on the closure of each complementary region of $|\lambda|$.

Roughly speaking, if $x, y ~(\in \h^2)$ are sufficiently close to each other, then the hyperbolic tangent planes of $\beta_L$ at $x$ and $y$ intersect at the external angle approximately equal to $\mu(\ol{x y})$ (with respect to the normal vector field of $\beta$; c.f. the hyperbolic tangent planes defined in \S \ref{ThP}).
Assume that $l$ is an isolated leaf of $\lambda$ and that $Q$ and $R$ are two adjacent complementary regions of $|\lambda|$ separated by $l$. 
Then $\beta_L(cl(Q))$ and $\beta_L(cl(R))$  are isometric copies of $cl(Q)$ and $cl(R)$ that intersect at the external angle $w(l)$.
This property determines the bending map $\beta_L$ if $L$ consists of isolated leaves.  
For general $L$, there is a sequence $(L_i)$ of measured laminations on $\h^2$ consisting of finitely many leaves, that approximates $L$ in Thurston's topology.
Then $\beta_{L} = \lim_{i \to \infty} \beta_{L_i}$ uniformly on compacts. 

Let $M$ be the sublamination of $L$ that consists of the leaves of $L$ with positive weight (note that $|M| \st \h^2$ has zero Lebesgue measure). 
Then, there is a unique tangent plane of $\beta_L$ at each point of $\h^2 \sm |M|$, and the tangent plan changes continuously (see \S \ref{ThP} and \cite{KP}). 
Therefore, the bending map $\beta_L$ is  $C^1$-smooth on $\h^2 \sm |M|$.

\subsection{Convex Hull Boundaries}\label{cb}(For details, see \cite[II.1.12 - 1.14]{EM}.)
Let $X$ be  a simply connected (open) region in $\RS$.
Then $X$ can be regarded as a projective structure on $\mr{\D}^2$. 
Consider $Conv(\RS \sm X)$, the convex hull of $\RS \sm X$.
It turns out that $\pt Conv(\RS \sm X)$ is isometric to $\h^2$ with respect to the induced path metric on $\pt Conv(\RS \sm X)$.
There is a unique measured lamination $L$ on $\h^2$ such that  $L$ does not contain leaves of weight more than $\pi$ and its bending map $\beta_L$ realizes the isometry from $\h^2$ to  $\pt Conv(\RS \sm X)$.
Then, by the orthogonal projection along geodesics in $\h^3$,  $~X$ maps onto $Im(\beta_L)$. 

\subsection{Thurston's Parameterization of Projective Structures}\label{ThP}

Let $P(S)$ be the space of all projective structures on $S$.  Thurston gave a
parametrization of $P(S)$ that reflects the geometry of projective structures in a combinatorial manner.
This parametrization is useful for the proof of the main theorems of this paper, since it involves a
decomposition of $\td{S}$ into $f$-injective regions, where $f: \td{S} \to \RS$ is the developing map of a projective structure. 

\begin{theorem}[Thurston] $P(S)$ is naturally homeomorphic to the product of  the Teichim\"uller Space of $S$ and the space of measured laminations on $S$:
 \begin{eqnarray} P(S) \simeq   \mathcal{T}(S) \times \mathcal{ML}(S) \ (\simeq
 \R^{6g-6} \times \R^{6g-6}) \label{eq1} \end{eqnarray}
\end{theorem}
\noindent (The proof is in \cite{KT}. For the following discussion, see also \cite{Ta}.)
Below we discuss some properties of this homeomorphism. 
An element in the left hand side of (\ref{eq1}) is a pair $(f, \rho)$, where $f: \td{S} \to \RS$ and $\rho: \pi_1(S) \to \PSL$.
Take a pair $(\tau, L) \in \mathcal{T}(S) \times \mathcal{ML}(S)$.
Then the total lift $(\h^2, \td{L})$ of $(\tau, L)$ induces a bending map $\beta: \h^2 \to \h^3$.
Furthermore, since the action of $\pi_1(S)$ preserves $\td{L}$, $~\beta$ induces a representation $\rho_L: \pi_1(S) \to \PSL$ such that $\beta$ is $\rho_L$-equivariant. 
If $(f, \rho)$ and $(\tau, L)$ represent the same projective structure, then $\rho = \rho_L$.
Letting $\rho_{id}: \pi_1(\oD) \to \PSL$ be the trivial representation, $(f, \rho_{id})$ is a projective structure on $\oD$, which is the universal cover of $C$.
Moreover, $(f, \rho_{id})$ corresponds to the measured lamination $(\h^2, \td{L})$ through the orthogonal projection and the bending map, which generalizes the correspondence between a simply connected region in $\RS$ and an injective bending map discussed in \S \ref{cb}.
Namely, in our current case, $f$ and $\beta$ are not necessarily embeddings, and we need to divide the domain of $f$ and the domain of $\beta$ so that their corresponding subdomains are homeomorphic through the orthogonal projection.

We shall first discuss the same correspondence for projective structures on $\oD$, which is more general than the above case.
Namely, there is a bijective correspondence between the projective structures on $\oD$ (that are not conformally equivalent to the Euclidian plane) and the measured laminations on $\h^2$ (up to the action of $\psl$); see \cite[Corollary 11.7]{KP}. 
For a measured lamination $L = (\lambda, \mu)$ on $\h^2$, let $C = C(L) = (f, \rho_{id})$ denote the corresponding projective structure on $\oD$. 
We shall discuss the correspondence between $C(L)$ and $L$.
There are a (topological) measured lamination $L' = (\lambda', \mu')$ on $(\oD, C)$ and the {\it collapsing map} $\kappa: (\oD, C, L') \to (\h^2, L)$, which describe the subdivision and the orthogonal projections. 
For each leaf $l$ of $L$ with positive weight, $~\kappa^{-1}(l)$ is a crescent of angle $w(l)$ with the canonical foliation (compare \cite[11.12]{K2}). 
Conversely, each crescent of angle $h$ in $(S, C)$ projects to a leaf of weight $h$ via $\kappa$ in the way discussed  in \S \ref{FlatCylinders}.
In the complement of such crescents, $\kappa$ is an isomorphism, i.e. a $C^1$-diffeomorphism that preserves the measured lamination. 
In summary, $L'$ is topologically obtained from $L$ by blowing up each leaf $l$ with positive weight of $L$ as above.
(Note that there is no periodic leaves with positive weight of $L'$.)
The collapsing map $\kappa$ is a continuous surjective map that homeomorphically takes each leaf of $L'$ to a leaf of $L$ and each component of $(\oD, C) \sm |L'|$ to a component of  $(\h^2, L) \sm |L|$.
Furthermore, this correspondence is bijective except the correspondence between the leaves of the crescents and the leaves of positive weight.

A {\it maximal ball} of a projective structure $C$ on $\oD$ is a maximal open subset of $\oD$ that $f$ homeomorphically takes to a round open disc in $\RS$, where the maximality is defined with respect to the set inclusion. 
If $U$ is a maximal ball, then $\pt f(U)$ is a round circle in $\RS$, and $Conv( \pt f(U)) \st \h^3$  is a copy of $\h^2$ whose ideal boundary is $\pt f(U)$. 
Let $H_U = Conv( \pt f(U))$, and let $\Psi_U: f(U) \to H_U$ be the orthogonal projection along geodesics in $\h^3$.  

Let $R$ be the closure of a component of  $\oD \sm |L'|$, or a leaf of  $L'$ that does not bound a component of  $\oD \sm |L'|$.
Then, $R$ is contained in a unique maximal ball $U$ and $R = R_U$ is called the {\it core} of $U$.
Conversely, each maximal ball $U$ contains a unique core.
These cores of maximal balls form a partition of $\oD$. 
Let $\beta: \h^2 \to \h^3$ be the bending map induced by $L$.
Then we have  $\Psi_U \circ f = \beta \circ \kappa$ on each core $R_U$, which describes the correspondence of $f$ and $\beta$. 
Define $\Psi: (\oD, C) \to \h^3$ by $\Psi(x) = \Psi_U(x)$ when $x \in R_U$. 
Then we have  $\Psi = \beta \circ \kappa$.

\begin{figure}[htbp]

\includegraphics[width = 4in]{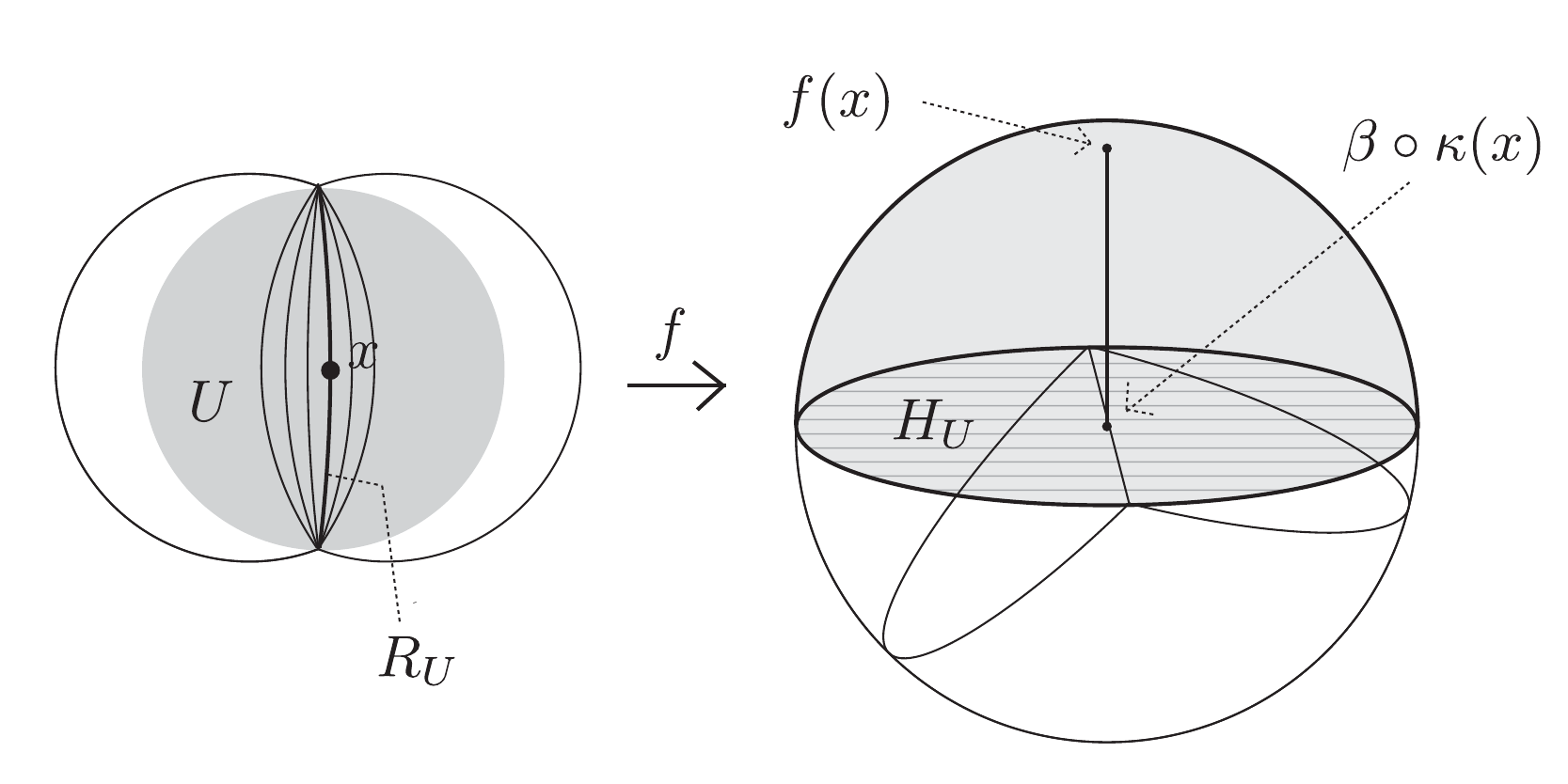}

\caption{}
\label{TangentPlane}
\end{figure}

Let $W$ be the union of leaves of $L$ with positive weight.  
Recall that $\beta$ is $C^1$-smooth except on $W$.
\begin{definition}
The {\bf hyperbolic tangent plane} of $\Psi$ at $x$ is $H_U = \pt Conv(f(U)) \cong \h^2$ when $x \in R_U$ (see Figure \ref{TangentPlane}).
\end{definition}
\noindent
This tangent plane is a support plane of $\beta(U_{\kp(x)})$ at $\Psi(x)$ where $U_{\kp(x)}$ is a sufficiently small neighborhood of $\kp(x)$.
Then this hyperbolic tangent plane  coincides with the standard hyperbolic tangent plane at each point of $\td{S} \sm \td{\kappa}^{-1}(W)$, which is the complement of the disjoint foliated crescents.
When $x \in \td{S}$ moves infinitesimally, the hyperbolic tangent plane of $\Psi$ at $x$ rotates about $\beta(l_x)$ in $\h^2$ by the amount of the transversal measure $\mu'$, where $l_x$ is a leaf of $L'$ through $x$ if it exists. 
In particular, when $x$ moves along a leaf or a moves in the closure of a component of $\td{S} \sm |L'|$, then the hyperbolic tangent plane does not change. 
Moreover, the hyperbolic tangent planes depend continuously on $x\in \tilde{S}$ (see \cite[Theorem 6.2]{KP}).

Let us return to the correspondence between a measured lamination $L$ on $(S, \tau)$ and a projective structure $C = (f, \rho)$ on $S$.
Since $\td{S} \cong \oD$, we have the canonical lamination $\td{L}'$ on $(\td{S}, \td{C})$ and the collapsing map $\td{\kp}: (\td{S}, \td{C}, \td{L}') \to (\h^2, \td{L})$.
By its construction, $\td{L}'$ is invariant under the action of $\pi_1(S)$, and it induces a measured lamination $L' = \td{L}'/ \pi_1(S)$ on $S$.
In addition, $\td{\kp}$ is $\rho$-equivariant, and it induces the {\it collapsing map} $\kappa: (S, C, L') \to (S, \tau, L)$.
Accordingly, for each periodic leaf $l$ of $L$, $~\kappa^{-1}(l)$ is a flat cylinder of height $w(l)$ with the canonical foliation. 
Conversely, each foliated flat cylinder of height $h$ in $(S, C)$ projects to a periodic leaf of weight $h$ via $\kappa$.
In the complement of such flat cylinders, $\kappa$ is an isomorphism. 

\begin{figure}[htbp]

\includegraphics[width= 5in]{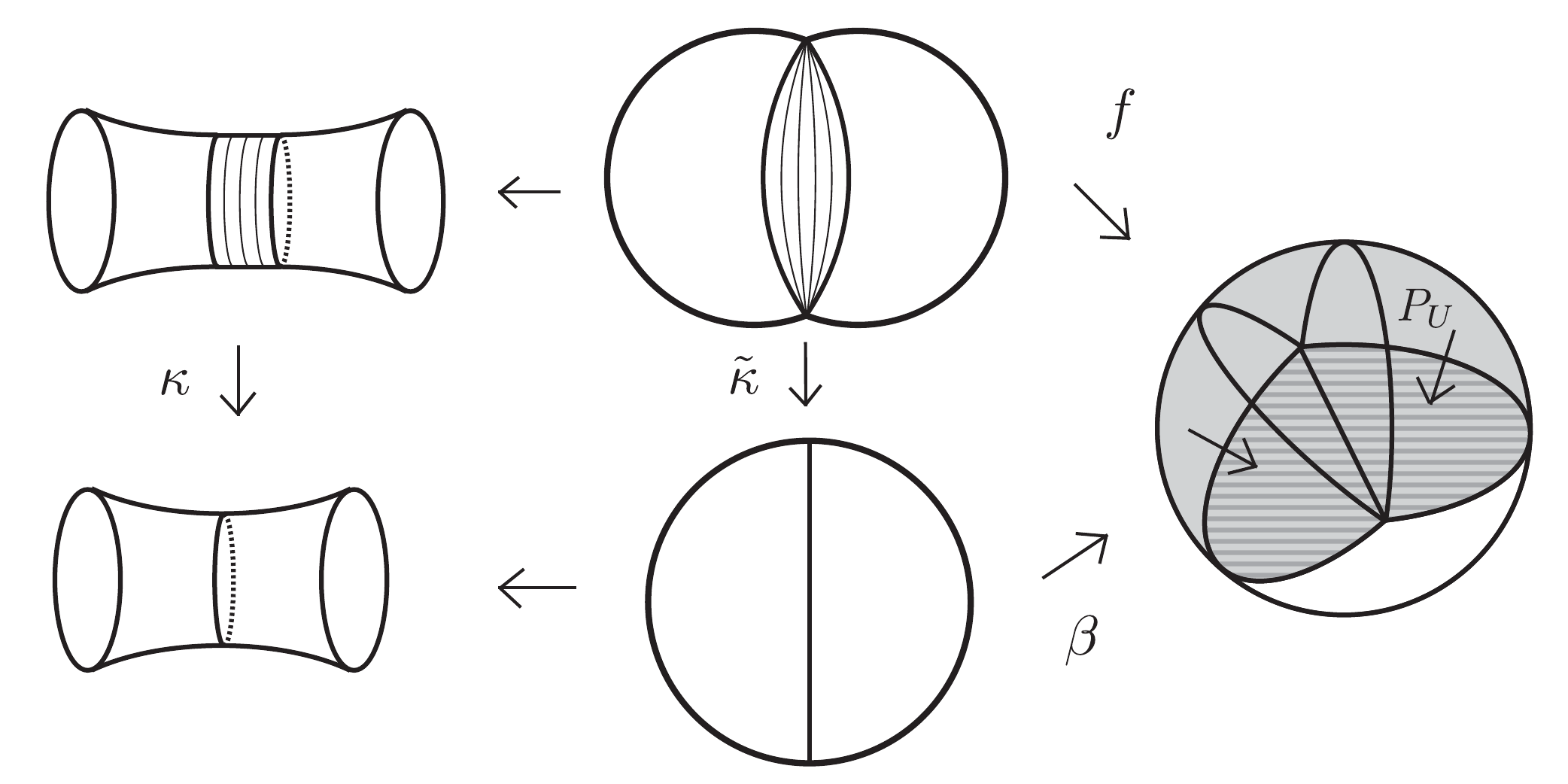}
\caption{}
\label{f&kappa&beta}
\end{figure}

Figure \ref{f&kappa&beta} illustrates the basic case when we have a measured lamination consisting of a periodic leaf on a complete hyperbolic cylinder. 
Near a periodic leaf of $L$ on $(S, \tau)$, we locally have a similar diagram.  

\subsubsection{The Intersection of a Lamination  and a Convex Set in $\h^2$}\label{intersection}

Let $L = (\lambda, \mu)$ be a measured lamination on $\h^2$.  
\begin{definition}
Let $X$ be a geodesic or a convex subset of $\h^2$ bounded by geodesics.
The {\bf intersection} of $L$ and $X$ is a measured lamination $(\lambda_X, \mu_X)$ on $\h^2$,  where $\lambda_X =\ cl\{ \ l \in \lambda \ | \ l \cap X \not= \emptyset \}$ and $\mu_X$ is the restriction of $\mu$ to $\lambda_X$. 
  Denote the intersection by $I(L, X)$.
\end{definition}

\begin{definition}
Let $X$ be a convex subset of  $\h^2$ bounded by geodesics.
The {\bf restriction} of $L$ to $X$ is  a measured lamination  $( \lambda|_X, \mu|_X )$ on $X$, where $\lambda_X = \{ l \cap X \>|\> l \in \lambda \}$ and $\mu|_X$ is defined by $\mu|_X(s) = \mu(s)$ for all geodesic segments $s$ in $X$.
Denote the restriction by $L|_X$.
\end{definition}

\begin{definition}
Let $C_1$ be a projective structure on a surface $F_1$, and let $F_2$ be a subsurface of $F_1$.
Then the {\bf restriction} of $C_1$ to $F_2$ is the restriction of the atlas of $C_1$ to $F_2$. 
Let $\mc{R}(C_1, F_2)$ denote the restriction of $C_1$ to $F_2$.
Conversely, if a projective structure $C_3$ on a surface $F_3$ is isomorphic to the one obtained by restricting $C_1$ to a subsurface of $F_1$, then we say that $C_3$ {\bf embeds} into $C_1$. 
\end{definition}

Let  $C_1 = (f_1, \rho_1)$ and assume, in addition, that the inclusion $F_2 \st F_1$ is $\pi_1$-injective.
Then the above definition is equivalent to the following: 
The restriction $\mc{R}(C_1, F_2)$ is the projective structure $(f_1|_{\td{F_2}}, \rho_1|_{\pi_1(F_2)})$, where $f_1|_{\td{F}_2}$  is the restriction of  $f_1$ to a lift $\td{F_2}$ of $F_2$ to  the universal cover of $F_1$ and  $\rho_1|_{\pi_1(F_2)}$ is the restriction of $\rho_1$ to $\pi_1(F_2)$ acting on $\td{F_2}$.

Let $L$ be a measured lamination on $\h^2$, and let $C(L) = (f_L, \rho_{id})$, the projective structure on $\oD$ corresponding to $L$ (\S \ref{ThP}). 
Let $X \,(\st \h^2)$ be either a geodesic or a convex subset bounded by geodesics. 
Let  $I = I(L, X)$, and let $C(I) = (f_I, \rho_{id})$.  
We also let $\kappa_L$ and $\kappa_I: \oD \to \h^2$ be the collapsing maps for $C(L)$ and $C(I)$, respectively. 

\begin{lemma}\label{C(I)/stC(L)}
There exists a homeomorphism $\phi:  \kp_I^{-1}(X) \to \kp_L^{-1}(X)$ such that $f_I = f_L \circ \phi$ on $\kp_I^{-1}(X)$.
Moreover, $C(I)$ embeds into $C(L)$.
\end{lemma}

\begin{proof}
Consider the leaves of $L$ and components of $\h^2 \sm |L|$ that intersect $X$, and let $X'$ be the union of these leaves and components. 
Then $X'$ is a convex subset of $\h^2$ containing $X$, and it is bounded by some leaves of $L$.
We also have $I = I(L, X) = I(L, X')$.
Therefore, it suffices to prove the lemma for $X'$.
Let $L'$ and $I'$ be the canonical measured laminations on $(\oD, C(L))$ and $(\oD, C(I))$, respectively.

Since  $~L|_{X'} = I|_{X'}$, we can assume that $~\beta_L = \beta_I$ on $X'$, where $\beta_L$ and $\beta_I$ are the bending maps induced by $L$ and $I$, respectively.
Therefore,  $\beta_L = \beta_I$ on $cl(X')$ by the continuity of bending maps.
Recall that $\kp_I^{-1}(X')$ and $\kp_L^{-1}(X')$ are obtained from $X'$ in the exactly same way, namely by blowing up the periodic leaves of  $L|_{X'} = I|_{X'}$.
Therefore, we have a canonical homeomorphism $\phi: cl( \kp_I^{-1}(X') ) \to cl( \kp_L^{-1}(X'))$ such that $\phi$ isomorphically takes $I' |_{\kp_I^{-1}(X')}$ to $L' |_ {\kp_L^{-1}(X')} $ and $\kappa_I = \kappa_L \circ \phi$ on  $cl(\kp_I^{-1}(X'))$.
Furthermore, the hyperbolic tangent plane of $\beta_I \circ \kp_I$ at $x \in cl(\kp_I^{-1}(X'))$ coincides with the hyperbolic tangent plane of $\beta_L \circ \kp_L$ at $\phi (x)$.
(See \S \ref{ThP}.)
The maximal ball of $C_I$ whose core contains $x$ maps to a round open disk by $f_L$, and the maximal ball of $C_L$ whose core contains $\phi(x)$ maps to a round open disk by $f_I$.
The convex hull boundaries of these open desks are the hyperbolic tangent planes of $\beta_L \circ \kp_L$ at $\phi(x)$ and of $\beta_I \circ \kp_I$ at $x$, and therefore, they must agree.
Then the round disks on $\RS$ also coincide. 
Recall that $f_I(x)$ and  $\beta_I \circ \kp_I(x)$ are connected by a geodesic in $\h^3$ orthogonal to the hyperbolic tangent plane, and so are $f_L(\phi(x))$ and $\beta_L \circ \kp_L (\phi (x))$.
Hence, $f_I(x) = f_L(\phi(x))$.

Each component $H$ of $\h^2 \sm X'$ is an open or closed half plane bounded by a leaf $l$ of $L$, and it does not contain leaves of $I$.
Then $~\kp_L^{-1}(H)$ is a component of $\oD \sm \kp_L^{-1}(X')$.  
This component is simply connected and bounded by a leaf $l_L$ of $L'$ that maps to $l$ via $\kp_L$.
Similarly, $\kp_I^{-1}(H)$ is a component of $\oD \sm \kp_I^{-1}(X')$.
This component is simply connected and bounded by a boundary curve $l_I$ of  $\kp_I^{-1}(X')$.
In addition, $l_I$ is a leaf of $I'$ or contained in a component of $\oD \sm |I'|$, and  it is homeomorphic to $l_L$ and $l$ via $\phi$ and $\kp_I$, respectively.
The leaf $l_L$ is contained in a unique maximal ball $U$ of $C(L)$, whose convex hull boundary is the hyperbolic tangent plane of $\beta_L \circ \kp_L$ at each point in $l_L$. 
Similarly, $l_I$ is contained in a unique maximal ball $V$ of $C(I)$, whose convex hull boundary is the hyperbolic tangent plane of $\beta_I \circ \kp_I$ at each point in $l_I$.
Thus these hyperbolic tangent planes are the same planes in $\h^3$.
Therefore, we can identify $V$ and $U$ by a $C^1$-diffeomorphism $\psi: V \to U$ such that $f_I = f_L \circ \psi$ on $V$.
Then $\psi = \phi$ on $\kp_I^{-1}(X') \cap V$.
Therefore, the embedding $\phi: \kp_I^{-1}(X) \to \kp_L^{-1}(X)  \st \oD$ extends to $\kp_I^{-1}(H)$, preserving $f_I = f_L \circ \phi$.
The different components of $\oD \sm \kp_I^{-1}(X')$ map to different components of $\oD \sm \kp_L^{-1}(X')$ by the extension. 
Hence, since $\kappa_I^{-1}(\h^2) = \oD$, we have an embedding of $C(I)$ into $C(L)$.
\end{proof}

Assume that $X$ has non-empty interior and no boundary geodesic of $X$ transversally intersects a leaf of $L$ with positive weight.
Then $cl(\kp_I^{-1}(X))$ and $ cl(\kp_L^{-1}(X))$ are  $C^1$-smooth subsurfaces of $\oD$. 
Thus we immediately obtain
\begin{corollary}\label{BasicCoro}
$\mc{R}(\,C(I), cl(\kp_I^{-1}(X))\, )$ and $\mc{R}(\,C(L), cl(\kp_L^{-1}(X)) \,)$ are isomorphic as projective structures.
\end{corollary}

\begin{figure}[htbp]
\includegraphics[width= 6in]{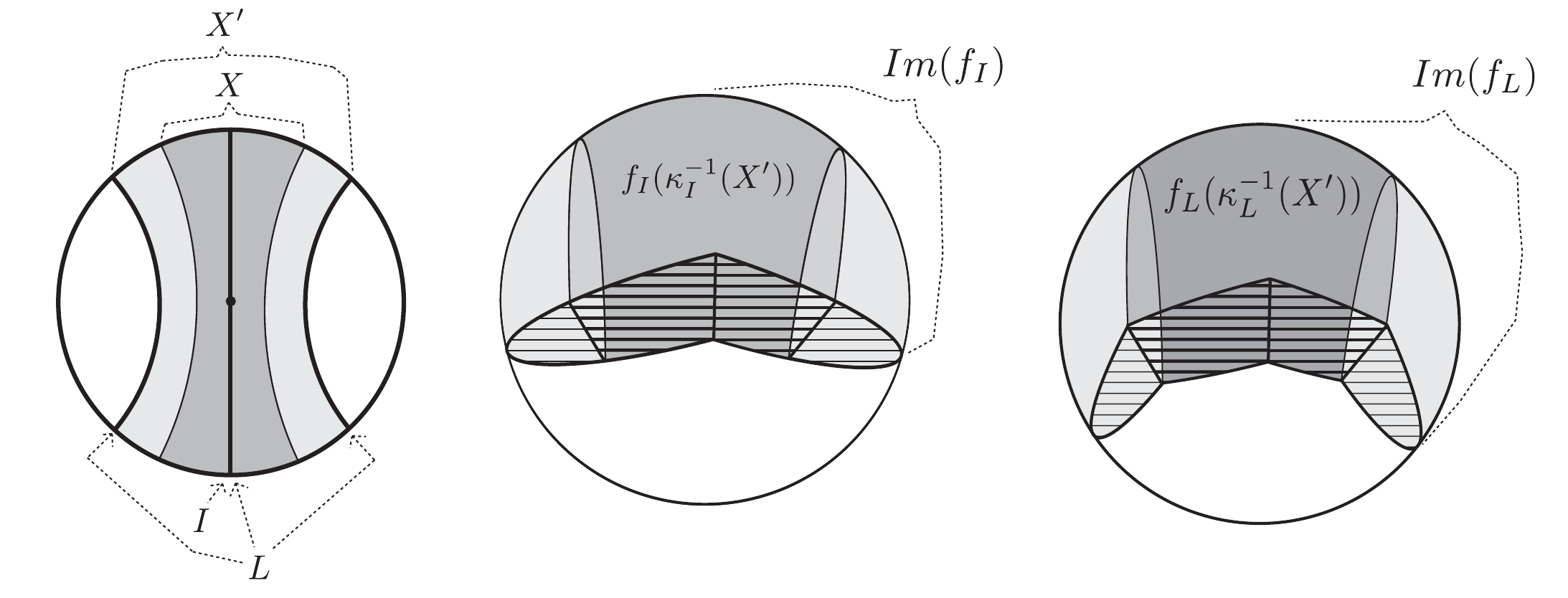}
\caption{A basic example for Lemma \ref{C(I)/stC(L)}.}
\label{BasicLemma}
\end{figure}

\subsection{A lemma on hyperbolic triangles}

Let  $\triangle ABC$ be an arbitrary hyperbolic triangle. 
Set  $a = length(BC)$, $b= length(CA)$, $c = length(AB)$, and  $\ap = \angle CAB$.
Then we have
\begin{lemma}\label{bilipschitz}
For every $\ep > 0$, there exists $S > 0$, 
such that $a > \frac{1}{2} (\sin \ep)\,(b+c)$ for every hyperbolic triangle $\triangle ABC$ with $\ap > \ep$.
\end{lemma}

\begin{proof}
In addition, we let $\bt = \angle ABC$,  and $\gm = \angle BCA$.

First, suppose that  $\ap \geq \pi/2$.
Then, since $\ap + \beta + \gm < \pi$, we have $\ap > \bt$ and $\ap > \gm$.
By the Hyperbolic Sine Rule, we see that $a > b$ and $a > c$.
Therefore $a > \frac{1}{2} (b + c) \geq \frac{1}{2} (\sin \ep) \, (b + c)$.

Next, suppose that  $0 < \ap < \pi/2$.
If $a \leq b$,
since $f(x) = \frac{\sinh x}{x}$ is a strictly increasing function, 
$$ \frac{\sinh a}{a} < \frac {\sinh b}{b}\,.$$
Then, using the Hyperbolic Sine Rule,
$$a > \frac{\sinh a}{\sinh b} \cdot b = \frac{\sin \ap}{\sin \bt} \cdot b \geq (\sin \ap) \cdot b\,.$$
This inequality $a >  (\sin \ap) \cdot b$ also holds when $a > b$. 
Similarly, we have $a > (\sin \ap) \cdot c\,$.
Since $0 <  \ap < \pi/2$, we have
$$a > \frac{1}{2} (\sin \ap) \, (b + c) > \frac{1}{2} (\sin \ep) \, (b + c)\,.$$
\end{proof}


\section{The Intersection of  a Lamination and Its Approximating Loop}\label{al}

Let $L = (\lambda, \mu)$ be an irrational minimal lamination on a closed orientable hyperbolic surface $S$, and $(l_i)$ be the sequence of simple loops on $S$ that converges to $|\lambda|$ constructed in \S \ref{Appro}. 
Let $\td{L} = (\td{\lambda}, \td{\mu})$ be the total lift of $L$ to $\h^2$, and let $\td{l}_i$ be a lift of $l_i$ to $\h^2$.
Let $L_i = (\lambda_i, \mu_i)$ be $I(\td{L}, \l_i)$ (see \S \ref{intersection}).
Note that the dual tree of $L_i$ is isometric to $\R$ 
(since $\td{l}$ intersects each leaf of $L_i$ in exactly one point).

\begin{definition}
Let $M = (\nu, \omega)$ be a measured lamination on a hyperbolic surface $F$.
Define the {\bf norm} of $M$ by 
$$ \| M \| = \sup\{ \omega(s)\},$$
where $s$ varies over all geodesic segments of length less than $1$ on $F$.  
\end{definition}
\noindent Remark: This norm is also called ``Thurston norm'' or ``roundness measure''.  

We next prove that the transversal measure $\mu_i$ of short geodesic segments is bounded by an arbitrary small number, provided that $i$ is large:
\begin{proposition}\label{propA}
$\displaystyle \lim_{i \to \infty} \| L_i \| = 0$.
\end{proposition}

The basic idea of the proof is that, when a geodesic segment $s$ with  $length(s) < 1$ intersects a fixed measured lamination at an angle close to zero, its transversal measure is also close to zero. 

Let $x$ be a point on a leaf $l$ of $\td{\lambda}$.
For $\theta \in ( -\pi/2, \pi/2]$,
let $l_{x,\theta}$ be the geodesic on $\td{S}$ intersecting $l$ at $x$ with $\angle_x(l, l_{x, \theta}) = \theta$ (see \S \ref{Appro}). 
Set $I(\tilde{L}, l_{x,\theta}) = (\lambda_{x,\theta}, \mu_{x,\theta})$.

 \begin{lemma}\label{ml2}
 For every $\epsilon > 0 $, there exists a constant  $\theta_0  > 0$ (which depends on $S, L, \epsilon$) such that, if $\theta \in ( -\theta_0, \theta_0)$ and $x \in | \td{\lambda} |$, then $\mu_{x, \theta}(s) < \epsilon$ for  all geodesic segments $s$ in $\h^2$ with $x \in s$ and $length(s) < 1$. 
\end{lemma}

\begin{figure}[htbp]
\includegraphics[width= 4in]{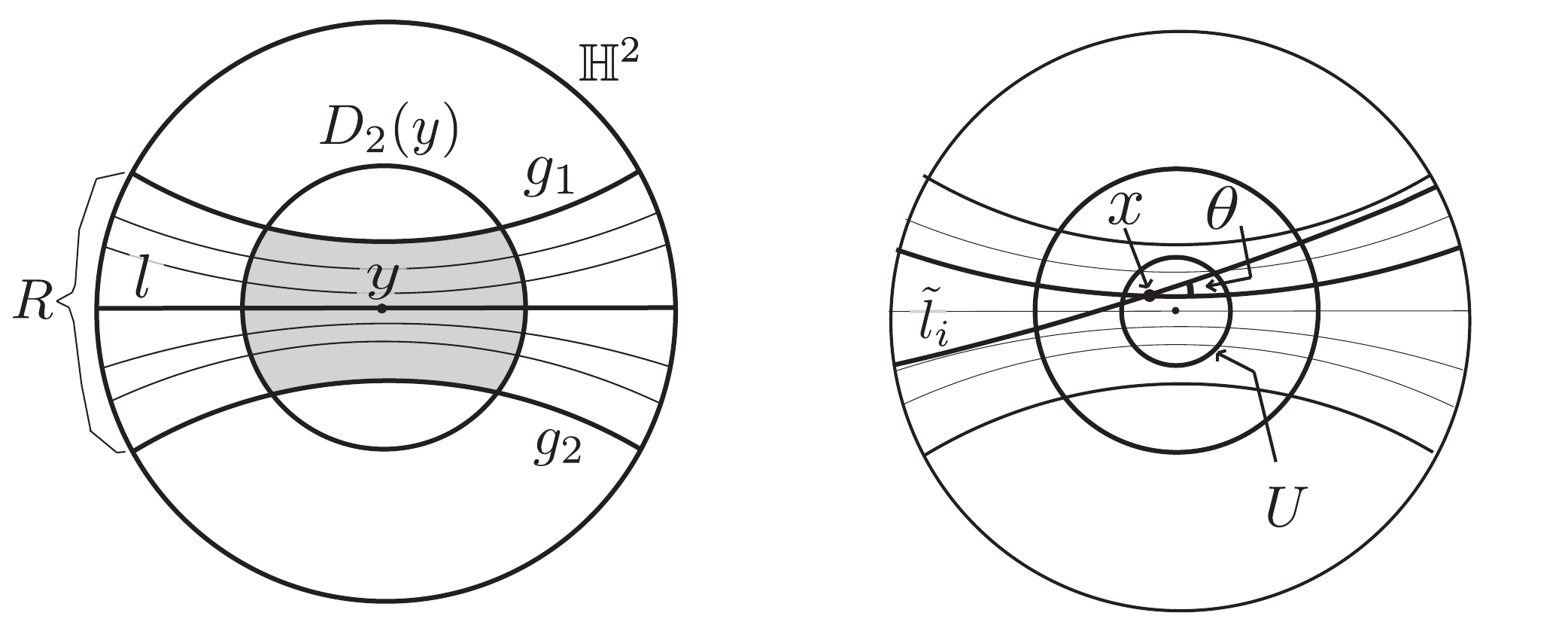}
\caption{In the left picture, $R \cap D_2(y)$ is shaded.}
\label{R-U-Theta}
\end{figure}

\begin{proof}
For an arbitrary  $y$ in $|\td{\lambda}|$, let $l$ be the leaf of  $\td{\lambda}$ through $y$. 
Then consider $D_2(y)$, the closed hyperbolic disk of radius $2$ centered at $y$.
In each component of $\h^2 \sm l$,  choose a geodesic $~g_i ~(i = 1, 2)~$ close to $l$ that does {\it not} transversally intersect a leaf of  $\td{\lambda}$,  i.e. $g_i$ is a leaf of $\td{\lambda}$ or is in the complement of $|\td{\lambda}|$.
Let $R \st \h^2$ be the convex region bounded by $g_1$ and $g_2$, which contains $l$.
For every $\epsilon > 0$, by applying Lemma \ref{box-nhd} to $l \cap D_2(y)$, we can assume that $g_1$ and $g_2$ are close enough to $l$, so that the $R \cap D_2(y)$ is contained in a flow box of height less than $\epsilon$.

Take a neighborhood $U$ of $y$ whose closure is contained in the interior of $R~ \cap ~D_1(y)$.
Then there exists  (small) $\theta_0 > 0$ such that, if $x \in U$ and $\theta \in (- \theta_0, \theta_0)$, then $l_{x, \theta} \st R$ (see  Figure \ref{R-U-Theta}).
Since $R \cap D_2(y)$ is contained in the flow box of height less than $\epsilon$ and $R$ supports $I(\td{L}, R)$,  for every geodesic segment $s$ in $D_2(y)$, the transversal measure of $s$ with respect to  $I(\td{L}, R)$ is bounded by $\epsilon$.
Therefore, since $I(\td{L}, l_{x, \theta})|_{D_2(y)}$  is a sublamination of $I(\td{L}, R)|_{D_2(y)}$, $~\mu_{x, \theta}(s) < \epsilon$.
If $s$ is a geodesic segment in $\h^2$ such that $s \cap U \neq \emptyset$ and $length(s) < 1$, then $s \st D_2(y)$. 
Thus, $~\mu_{x, \theta}(s) < \epsilon$.
This proves the lemma if $x$ is in $U$, which is a neighborhood of $y$.
Since $S$ is compact and $\td{L}$ is invariant under the deck transformations, the lemma follows.  
\end{proof}

\begin{figure}[htbp]
\includegraphics[width= 2.5in]{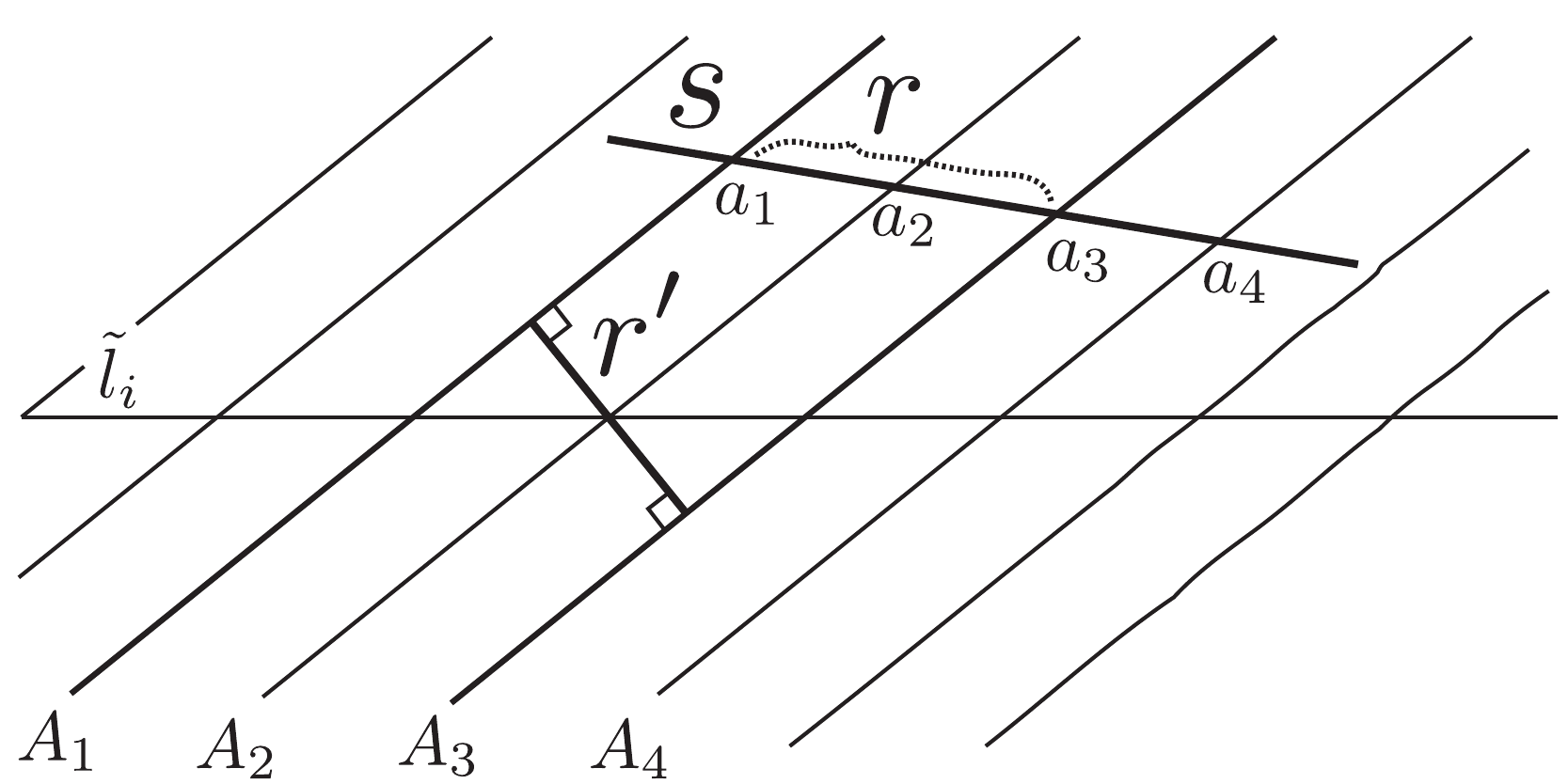}
\caption{A case that $p = 4$ and $q = 1$.}
\label{s&r}
\end{figure}

\begin{proof}(Proposition \ref{propA})
Fix arbitrary $\epsilon > 0$.
It suffices to show that, for sufficiently large $i$, if a geodesic segment $s$ in $\h^2$ satisfies $\mu_i(s) \ge \epsilon$, then  $length(s) \geq 1$. 
By Lemma \ref{ml2}, there exists $\theta_0 > 0$,   such that,  if $x \in |\td{\lambda}|$ and $\theta \in (-\theta_0, \theta_0)$, then $\mu_{x,\theta}(s) < \frac{2}{5}\epsilon$  for every geodesic segment $s$ with $x \in s$ and  $length(s) < 1$.
By Lemma \ref{ml1}, for sufficiently large  $i$, $~\angle(\td{\lambda}, \td{l}_i) = \angle (\lambda, l_i) < \theta_0$ and $\mu(l_i) < \epsilon/2$.

Consider the cyclic subgroup of  ${\rm PSL}({\rm 2}, \R)$ generated by the translation along $\td{l}_i$ by $length(l_i)$, which we can regard as $\pi_1(l_i) ( \st \pi_1(S) )$ acting on $\h^2$.
Choose a leaf $m \in \lambda_i$, and consider the orbit of $m$ under the action of the cyclic group.
This orbit forms a sublamination $\nu$ of $\lambda_i$.
The leaves of $\nu$ intersect $\td{l}_i$ at a constant angle less than $\theta_0$, and the distance between their consecutive intersection points is equal to $length(l_i)$. 
Since $\mu(l_i) < \epsilon/2$ and the dual tree of $L_i$ is isometric to $\R$, if a geodesic segment $s$ lies between two adjacent leaves of $\nu$,  then $\mu_i(s) < \epsilon/2$. 
Therefore, if $\mu_i(s) \ge \epsilon$, then $s$ transversally  intersects $|\nu|$ at least three times. 
Let $a_1, a_2, \ldots, a_p$ be the intersection points lying on $s$ in this order, and let $A_1, A_2, \ldots, A_p$ be the  leaves of $\nu$ through these points.
Take $q \in \N$ such that $2q +1$ is the maximal odd integer {\it not} exceeding $p$.
Let $r$ be the subsegment of $s$ with end points $a_1$ and $a_{2q+1}$.
(See Figure \ref{s&r}.)
Then $length(r) \leq length(s)$ and $\frac{2}{5} \mu_i(s) \leq  \mu_i(r) \leq \mu_i(s)$.
Let $r'$ be the geodesic segment that realizes the distance between $A_1$ and $A_{2q+1}$.
Then $r'$ is orthogonal to $A_1$ and $A_{2q+1}$.
In addition, $r'$ intersects $\tli$ transversally, since otherwise $\tli , r', A_1, A_{2q+1}$ bound a hyperbolic triangle whose interior angle sum is $\pi$ or a hyperbolic rectangle whose interior angle sum is $4\pi$. 
Therefore,  the triangle bounded by $\tli$, $r'$, $A_1$ is isometric to the triangle bounded by $\tli, r', A_{2q + 1}$.
Thus  $A_{q + 1} \cap \tli$ is the middle point of $r'$.
Note that $A_{q + 1} \cap \tli \in r'$, $~\tli$ intersects $|\nu| \st |\td{\lambda}|$ at  $A_{q + 1} \cap \tli$ at an angle less than $\theta_0$, and $\mu_i(\gamma') = \mu_i(\gamma) \geq \frac{2}{5}\epsilon$.
Therefore, $length(r') \geq 1$.
Hence, $1\leq length(r') \leq length(r) \leq length(s)$. 
\end{proof}

\section{Injectivity of Bending Maps}\label{IBM}

In this section, let $L = (\lambda, \mu)$ be a measured lamination on $\h^2$ and $\beta_L = \beta: \h^2 \to \h^3$ be the bending map induced by $L$.
Recall that $\|L\| = sup \{\ {\mu}(s)\  \}$,  where $s$ varies over all geodesic segments on $\h^2$ of length less than 1. 

\begin{theorem}[Epstein, Marden and Markovic; ~\cite{EMM}, Theorem 4.2.2]\label{ib} There exists $\delta \in (0, \pi)$ such that, if $\|L\| 
< \delta$, then the induced bending map $\beta_L$ is a bilipschitz embedding; hence, it continuously extends to $\partial \h^2$ as an embedding whose image is a simple loop in $\partial \h^3$.  \end{theorem}

To prove the decomposition theorem (Theorem \ref{id}), we need  a generalization of Theorem \ref{ib}. 
However, if we do not require flat cylinders in Theorem \ref{id} to be integral, then Theorem \ref{ib} is sufficient.  
To state the generalization, let us set up some notions. 
A leaf $l$ of $L$ is {\it outermost}, if the other leaves of $L$ lie only in one component of $\h^2 \sm l$.
Consider all outermost leaves of $L$ with positive weight.
Then this set forms a sublamination $\pt \lambda$ of $\lambda$, and $\pt \lambda$ consists of isolated leaves. 
Let $\pt L$ be the measured lamination on $\h^2$ obtained by assigning each leaf of $\pt \lambda$ its positive weight with respect to $\mu$. 
Let $L' = (\lambda', \mu')$ be a sublamination of $\pt L$. 
We also let $d (L') =\inf\{ dist_{\h^2}(l_1, l_2) ~|~ l_1, l_2 \in \lambda' , l_1 \neq l_2\} \geq 0$.
Let $L \sm  L' = ( cl(\lambda \sm \lambda'), \mu - \mu')$.   
Note that $|L'|$ bounds a convex region of $\h^2$ and that the convex region contains $|L \sm L'|$.

\begin{theorem}\label{ib2} 
 For every $D > 0$, there exists $\delta \in (0, \pi)$ with the following property:
If a measured lamination $L$ on $\h^2$ contains a sublamination $L'$ such that\\
(i) $L' \st \pt L $, \\
(ii)  $\| L \sm L' \| < \delta$,\\
(iii) $d(L') > D$, and\\
(iv) every leaf of $L'$  has weight less than $\pi/2$,\\
then the induced bending map $\beta_L\cn \h^2 \to \h^2$ is a bilipschitz map. 
\end{theorem}

\begin{proof}
Since $dist(x,y) \geq dist(\beta_L(x), \beta_L(y))$ for all $x, y \in \h^2$, it suffices to show that there exists $S > 0$ such that $dist(\beta_L(x), \beta_L(y)) > S \cdot dist(x, y)$ for all distinct $x, y \in \h^2$.
Let $p\cn [0, P] \to \h^2$ be the geodesic segment connecting $x$ to $y$, parametrized by arc length, where $P = dist(x, y)$. 
Let $P_0, P_1, \ldots, P_n$ denote the distinct points lying on $p$ in this listed order, such that $P_0 = x$, $P_n = y$, and $P_1, P_2, \ldots P_{n-1}$  are the transversal intersection points of $p$ and $|L|$. 
Let $R$ be the closed convex region bounded by $L'$.
Note that $\beta_L = \beta_{L \sm L'}$ on $R$.
Then, the following lemma immediately follows from Assumption (ii) (see Lemma 4.4 and the proof of Corollary 4.5 in \cite{EMM}):  \begin{lemma}\label{bounds}
For every $\epsilon > 0$, there exists $\delta > 0$ (depending only on $\epsilon$) such that
\begin{eqnarray*}
dist(\beta_L(x), \beta_L(y)) > (1 - \ep) ~ dist(x, y) ~ and\\
\angle \beta_L(x) \beta_L(y)\beta_L(P_{n-1}) < \ep 
\end{eqnarray*}
 for all distinct $x, y \in R$ (see Figure \ref{Angle&Length}).
\end{lemma}

\begin{figure}[htbp]
\includegraphics[width= 2.5in]{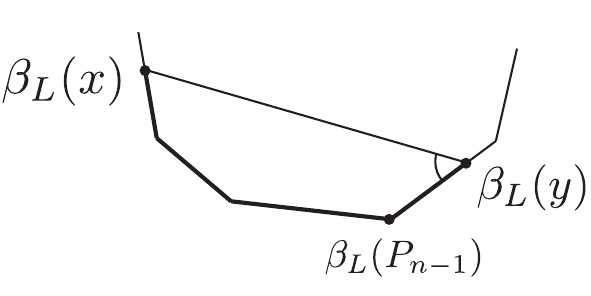}
\caption{}
\label{Angle&Length}
\end{figure}

{\bf Discrete Case.} First, we assume that $L$ contains only isolated leaves. 
Fix $\delta > 0$, obtained by applying Lemma \ref{bounds} to an arbitrarily fixed $\epsilon < 1$. 

\emph{Case 1.}   
Suppose that $x$, $y \in R$. 
Then, by Lemma \ref{bounds}, we have $dist(\beta_L(x), \beta_L(y)) > \epsilon ~ dist(x,y)$.

\emph{Case 2.}\ (See Figure \ref{Case2&3}.)
Suppose that $x \in R$ and $y \in \h^2 \sm R$.
Then, since $x, P_{n-1} \in R$, by Lemma \ref{bounds}, we have  $dist(\beta_L(x), \beta_L(P_{n-1})) > (1 - \ep) \, dist(x, P_{n-1})$ and $\angle \beta_L(x) \beta_L(P_{n-1}) \beta_L(P_{n-2}) < \epsilon$. 
Since $P_{n-1}$ is contained in a leaf $l$ of $L'$, by Assumption (iv) and the Triangle Inequality, we have $$\angle \bt_L(P_{n-2}) \bt_L(P_{n-1}) \bt_L(y) \geq \pi - w(l) >  \pi/2.$$
Therefore $$\angle \beta_L(x) \beta_L(P_{n-1}) \beta_L(y) > \pi - (\ep + \pi/2) = \pi/2 - \ep\,.$$

Thus, by Lemma \ref{bilipschitz},  there exists $S' > 0$, which does $not$ depend on the choices of $L, x, y$ (under the given conditions), such that
$$dist(\bt_L(x), \bt_L(y)) > S' (\,dist(\beta_L(x), \beta_L(P_{n-1}))\, + \,dist(\bt_L(P_{n-1}), \bt_L(y))\,).$$
Therefore $dist(\bt_L(x), \bt_L(y)) > S' (1 - \ep) dist(x,y).$

\emph{Case 3.}\ (See Figure \ref{Case2&3}.)
Suppose that $x, y \in \h^2 \sm R$.
Since $P_1 \in R$ and $y \in \h^2 \sm R$, by Case 2, there exists $S' > 0$, which does $not$ depend on the choices of $L, x, y$, such that
\begin{equation}
dist(\bt_L(P_1), \bt_L(y)) > S' \, dist(P_1, y). \label{Case3Eqn}
\end{equation}
By an argument similar to that in Case 2, we have 
$$dist(\beta_L(P_1), \bt_L(P_{n-1})) > (1 - \ep)\, dist(P_1, P_{n-1}) > (1 - \ep) D$$
and $\beta_L(P_1) \beta_L(P_{n-1}) \beta_L(y) > \pi/2 - \ep$.
Therefore, by taking a smaller $\epsilon > 0$ if necessary, we can assume that $\angle \bt_L(y) \bt_L(P_1) \bt_L(P_{n-1}) < \pi/2 - 2\ep$.
Thus, by Assumption (iv), we have $\angle \beta_L(x) \beta_L(P_1) \beta_L(y) > \pi - ((\pi/2 - 2\ep) + \pi/2 + \ep) > \ep$.
Then, by applying Lemma \ref{bilipschitz} to $\triangle \beta_L(x)\beta_L(P_1)\beta_L(y)$, we have $S'' > 0$ such that
$$dist(\bt_L(x), \bt_L(y)) > S''\, [\,dist(\beta_L(x),\beta_L(P_1)) + dist(\beta_L(P_1), \beta_L(y))\,].$$

Combining this inequality with (\ref{Case3Eqn}), we obtain $S > 0$ such that $dist( \beta_L(x),\beta_L(y)) > S \, dist(x, y)$
for all distinct $x, y \in \h^2 \sm R.$

\begin{figure}[hbt]
\includegraphics[width= 5in]{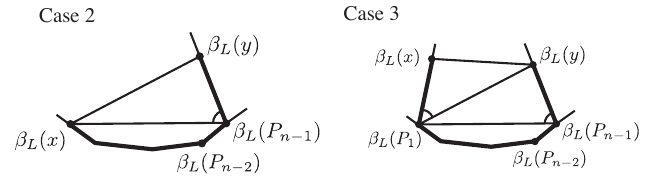}
\caption{} 
\label{Case2&3}
\end{figure}

{\bf General Case.}
Assume that $L$ does {\it not} consist of isolated leaves.
Then there is a sequence of measured laminations on $\h^2$, $(L_i = (\lambda_i, \mu_i))$, limiting to $L$ in Thurston's topology, such that each $L_i$ consists of finitely many leaves,  contains a sublamination $L_i' \st \pt L_i$, and  satisfies Assumptions (i) and  (iv).
Let $\beta_{L_i}: \h^2 \to \h^3$ denote the bending map induced by $L_i$. 
Then $(\beta_{L_i})$ converges to $\beta_L$ uniformly on compacts. 
(See \cite[II.1.13]{EM} for the construction of $(L_i)$ and the convergence of $(\beta_{L_i})$.)
By the Discrete Case, there exists $S > 0$, such that if $x, y \in \h^2$, then  
$$dist(\beta_{L_i}(x), \beta_{L_i}(y)) > S \, dist(x, y)$$
 for sufficiently large $i = i(x,y)$.
In addition, the convergence of $(\beta_{L_i})$ to $\beta_L$ implies that $(\,dist(\beta_{L_i}(x), \beta_{L_i}(y))\,)$ limits to $dist(\beta_L(x), \beta_L(y))$ as $i \to \infty$.
Then 
$$ dist(\beta_L(x), \beta_L(y)) \geq S \, dist(x, y)$$
for all $x, y \in \h^2$. 
\end{proof}

Assume that $L$ satisfies the assumptions of Theorem \ref{ib} or Theorem \ref{ib2}. 
By these theorems, $~\beta = \beta_L: \h^2 \to \h^3$ is an injective quasiisometric embedding, and, hence, it extends continuously to a homeomorphism $\pt \beta$  from  $\pt \h^2$ onto a simple loop on $\pt \h^3 \cong \RS$ (see \cite{Ghys}, \cite{Gromov}). 
Therefore, $\RS \sm Im(\pt\beta)$ consists of two simply connected regions.   

\begin{corollary} \label{admissibleD}

Under the assumption of Theorem \ref{ib} or Theorem \ref{ib2},
the projective structure $C(L)$ on $\mb{D}^2$ corresponding to $L$ is admissible. 
\end{corollary}

\begin{proof}
Since $\beta$ is an injective continuous quasiisometric embedding, $Im(\beta)$ is a proper surface embedded in $\h^3$.  
Therefore, $Im(\beta)$ separates $\h^3$ into two components (the Jordan-Brouwer Separation Theorem). 
Since $Im(\beta)$ is locally convex,  one of the components of $\h^3 \sm Im(\beta)$ is convex (see \cite[I.1.3]{CEG}). 

The concave component of $\h^3 \sm Im(\beta)$ is cobounded by $Im(\beta)$ and a topological closed disk $D$ contained in $\RS$.
Then $Im(\beta) = \pt Conv(\RS \sm D)$.
Since $L$ does not contain leaves with weight $\geq \pi$, $~L$ is the canonical bending lamination on $\pt Conv(\RS \sm D)$.
Therefore, $(\h^2, L)$ is Thurston's coordinates for the projective structure on $\mr{D}$; Hence $C(L)$ is admissible.
\end{proof}
 

\section{The Existence of Admissible Loops}

Let $S$ be a closed orientable surface of genus at least $2$.
Let $C = (f, \rho)$ be a projective structure on $S$.
Express $C$ in Thurston's  coordinates as $(\tau, L)$, where $\tau$ is a marked hyperbolic structure on $S$ and $L = (\lambda, \mu)$ is a measured lamination on $(S, \tau)$.
Let us recall other related notions from \S \ref{P}:
Let $\td{L} = (\td{\lambda}, \td{\mu})$ be the total lift of $L$ to $\h^2$, and
let $\beta_{\td{L}} : \h^2 \to \h^3$ be the bending map induced by $\td{L}$.
Let $\kappa: (S, C) \to (S, \tau)$ be the collapsing map, and let $\td{\kappa}: (\td{S}, \td{C}) \to (\h^2, \td{L})$ be the lift of $\kp$ to a map between the universal covers of $(S, C)$ and $(S, \tau)$.
Let $L' = (\lambda', \mu')$ be the canonical (topological) measured lamination on $(S, C)$ corresponding to $L$ via $\kappa$ (\S \ref{ThP}). 
Let $\td{L}'$ be the total lift of $L'$ to $(\td{S}, \td{C})$.

Since $L$ decomposes into minimal measured laminations, we can set
$$L = ( P_1 \,\bigsqcup\,~ P_2 \,\bigsqcup \,  \ldots \,\bigsqcup \, P_m ) \,\bigsqcup \, ( M_1 \,\bigsqcup \, M_2 \,\bigsqcup \, \ldots \,\bigsqcup \, M_n),$$
where $P_h, ~ h=1, \ldots, m,$ are the periodic minimal sublaminations of $L$ and  $M_j = (\nu_j, \omega_j), ~j= 1, \ldots, n,$ are the irrational minimal sublaminations of $L$.
Let $p_h = |P_h|$ denote the periodic leaf supporting $P_h$, 
 and let $p$ denote the periodic part of $|\lambda|$, 
 $$p_1 \,\bigsqcup \, p_2 \,\bigsqcup \, \ldots \,\bigsqcup \, p_m.$$
Let $M = (\nu, \omega)$ denote  the irrational part of $L$, 
$$ M_1 \,\bigsqcup \, M_2 \,\bigsqcup \, \ldots \,\bigsqcup \, M_n.$$
For each $j \in \{1,2, \ldots, n\}$, let $(l_{i,j})^\infty_{i =1}$ be the sequence of simple geodesic loops on $(S, \tau )$ that approximates $|\nu_j|$, constructed in \S \ref{Appro}. 
By Lemma \ref{ml1}, we have $\lim_{i \to \infty} \omega_j(l_{i,j}) = 0$.
Since $|\nu_j|$ is an isolated subset of $|\lambda|$, we can assume that $l_{i,j}$ does not intersect $|\lambda \sm \nu_j|$.
Therefore, $\lim_{i \to \infty} \mu(l_{i,j}) = 0$.
Let $l_i = l_{i,1}\,\bigsqcup \, \ldots \,\bigsqcup \, l_{i,n}\,$, so that $\lim_{i \to \infty} l_i = |\nu|$.

Recall that, for each $h$, $~\kappa^{-1}(p_h)$ is a flat foliated cylinder of height $w(p_h)$ in $(S, C, L')$, where $w(p_h)$ is the weight of $p_h$. 
The foliation on $\kappa^{-1}(p_h)$ consists of admissible loops that are homeomorphic to $p_h$ via $\kappa$ (\S \ref{ThP}).

Recall also that $\kappa$ is a $C^1$-diffeomorphism on $S \sm\kappa^{-1}(p)$.
Therefore, if  $l$ is an essential simple loop on $(S, \tau)$ disjoint from $|\lambda|$,  then $\kappa^{-1}(l)$ is an essential simple loop on $(S, C)$ disjoint from $|\lambda'|$.
We shall see that $\kappa^{-1}(l)$ is also admissible. 
Let $P$ be the component of $(S, \tau) \sm |L|$ containing $l$, and let $\td{l}$ and $\td{P}$ be the lifts of $l$ and $P$, respectively.
We can assume that  $\td{l} \st \td{P}$ and that $\pi_1(l)$ acts on $\td{l}$ freely and properly discontinuously as an infinite cyclic subgroup of $\PSL$ generated by a hyperbolic element.
Let $\gamma_l$ be the homotopy class of $l$ that generates $\pi_1(l)$.
Since $\td{l}$ is a quasigeodesic in $\h^2$ and $\beta_{\td{L}}$ is an isometry on $\td{P}$,  $~\beta_{\td{L}}(\td{l})$ is a quasigeodesic in $\h^3$.
Then $\pi_1(l)$ acts on $\beta_{\td{L}}(\td{l})$ freely and properly discontinuously via $\rho$, and this action fixes the distinct end points of $\beta_{\td{L}}(\td{l})$ on $\RS$.
Therefore, $\rho(\gamma_l)$ is a loxodromic element of $\PSL$. 
The curve $\td{\kappa}^{-1}(\td{l})$ is a lift of $\kp^{-1}(l)$ to $(\td{S}, \td{C})$.
Then $\td{\kappa}^{-1}(\td{l})$ is contained in $\td{\kp}^{-1}(\td{P})$, which is a component of $(\td{S}, \td{C}) \sm |\td{L}'|$.
Since $cl(\td{\kp}^{-1}(\td{P}))$ is the core of a maximal ball,  $\td{\kappa}^{-1}(\td{l})$ is contained in a maximal ball. 
Thus, $f$ is an embedding on $\td{\kappa}^{-1}(\td{l})$.
Hence, $\kappa^{-1}(l)$ is admissible. 

We have  $\mu(p_h) = 0$ and $\mu(l) = 0$.
For each $i, j$, let $l'_{i,j} = \kappa^{-1}(l_{i,j})$, which is a simple loop on $(S, C)$.
The fact that $\lim_{i \to \infty} \mu(l_{i,j}) = 0$ suggests the following proposition:

\begin{proposition}
For each $j \in \{ 1, 2, \ldots, n \}$, $~l'_{i,j}$ is admissible, provided that $i$ is sufficiently large.
\end{proposition}

\begin{proof}
Let $\tl'_{i,j}$ be a lift of $l_{i,j}$ to $\h^2$.
Consider the measured lamination $I(\td{L}, \tl'_{i,j})$.
Since $l_{i,j}$ is disjoint from $|\lambda \sm \nu_j|$,  $~I(\td{L}, \tl_{i,j}) = I(\td{M}_j, \tl'_{i,j}) $, where $\td{M}_j$ is the total lift of $M_j$.
Choose $\delta > 0$ as in Theorem \ref{ib}. 
Applying Proposition \ref{propA} with $I(\td{L}, \tl'_{i,j}) = I(\td{M}_j, \tl'_{i,j})$, we have $~ \| I(\tilde{L}, \tl'_{i,j}) \| < \delta$ for all large $i$.
By Theorem \ref{ib}, for sufficiently large $i$, the bending map $~\beta_{I(\tilde{L}, \tl'_{i,j})}$ induced by $I(\tilde{L}, \tl'_{i,j})$ is an injective quasiisometric embedding, and it continuously extends to an embedding of $\pt \h^2$.

Let  $\gamma_{i,j}$ be the homotopy class of $l_{i,j}$ that acts on $\h^2$ as a hyperbolic element of $\psl$  preserving $\tl'_{i,j}$.
The extension of $\beta_{I(\tilde{L}, \tl'_{i,j})}$ homeomorphically takes the limit set of $\langle \gamma_{i,j} \rangle$ to the limit set of $\langle \rho(\gamma_{i,j}) \rangle$.
Thus $~\rho(\gamma_{i,j})$ is loxodromic. 

By Corollary \ref{admissibleD}, the projective structure $C(I(\td{L}, \tl'_{i,j}))$ on $\oD$ corresponding to $I(\td{L}, \tl'_{i,j})$ is admissible  for sufficiently large $i$. 
Then the developing map $f_I$ of $C(I(\td{L}, \tl'_{i,j}))$ is an embedding.
By Lemma \ref{C(I)/stC(L)}, there exists a homeomorphism $\phi:  \kappa_I^{-1}(\tl_{i,j}) \to \kp^{-1}(\tl_{i,j})$ such that $f_I = f \circ \phi$ on $\kappa_I^{-1}(\tl_{i,j}) $, where $\kp_I: \oD \to \h^2$ is the collapsing map for  $C(I(\td{L}, \tl'_{i,j}))$. 
Since $f_I$ is an embedding, $f = f_I \circ \phi^{-1}$ restricted to $\tl'_{i,j}$ is an embedding.  
\end{proof}

We thus obtain an admissible loop from every minimal sublamination of $L$ and every complementary region of $|\lambda|$ that is not topologically an open disk.  Therefore, 

\begin{theorem}\label{gc}
For every projective structure $C$ on $S$, there exists an admissible loop on $(S, C)$.
 \end{theorem}

\noindent {\it Remark:} Equivalently, we can state that every projective structure on $S$ admits a grafting operation (see \cite{Go2} for the definition of a grafting  operation). 


\section{Admissible Decomposition}
We carry over our notation from the previous section.
We have shown that $l^j_i$ and $p_h$ correspond to admissible loops on $S$ through $\kappa$, provided that $i$ is sufficiently large.
Their union 
$$l_i \,\bigsqcup \, p  = (l_{i,1} \,\bigsqcup \, \ldots \,\bigsqcup \, l_{i,n}) \,\bigsqcup \,  (p_1\,\bigsqcup \, \ldots \,\bigsqcup \, p_m)$$
 is a multi-loop on $(S, \tau)$.
In this section, we show that $l_i \bigsqcup p$ decomposes $(S, C)$ into admissible subsurfaces.  

\begin{theorem}[Admissible Decomposition]\label{id}
Let $C$ be a projective structure on a closed orientable surface $S$ of genus at least 2. Then there exists a decomposition of $S$ into  cylinders and compact subsurfaces of negative Euler characteristic, such that the restriction of $C$ to each cylinder is an integral flat structure and the restriction to each subsurface of negative Euler characteristic is an admissible projective structure.  
\end{theorem}   

Note that every flat cylinder of height less than $2\pi$ is admissible. 
Therefore, every integral flat cylinder can be further decomposed into admissible flat cylinders, if we wish. 
Moreover, by further decomposing each surface of negative Euler characteristic into pairs of pants, if necessary, we immediately obtain:  

\begin{corollary}\label{Pants}
There exists a decomposition of $S$ into pairs of pants and cylinders such that the restriction of $C$ to each cylinder is an integral flat structure and the restriction to each pair of pants is an admissible structure. 
\end{corollary}

Let $l$ be a geodesic lamination on a complete hyperbolic surface $F$. 
Let $\mc{NT}(l)$ denote the collection of all geodesic segments of length less than one on $F$ that do {\it not} transversally intersect any leaves of $l$.  
Then,  a geodesic segment $s$ connecting $x$ and $y$ on $F$ is an element of $\mc{NT}(l)$ if and only if either $s \st |l|$ or $(s \sm \{x, y\}) \cap |l| = \emptyset$.

\begin{lemma}\label{LemId}
For every $\epsilon > 0$, there exists $i_0 \in \N$ such that, if  $i > i_0$, then $\omega(s) < \epsilon$ for all $s \in  \mathcal{NT}(l_i)$. 
\end{lemma}

\begin{figure}[htbp]
\includegraphics[width= 4in]{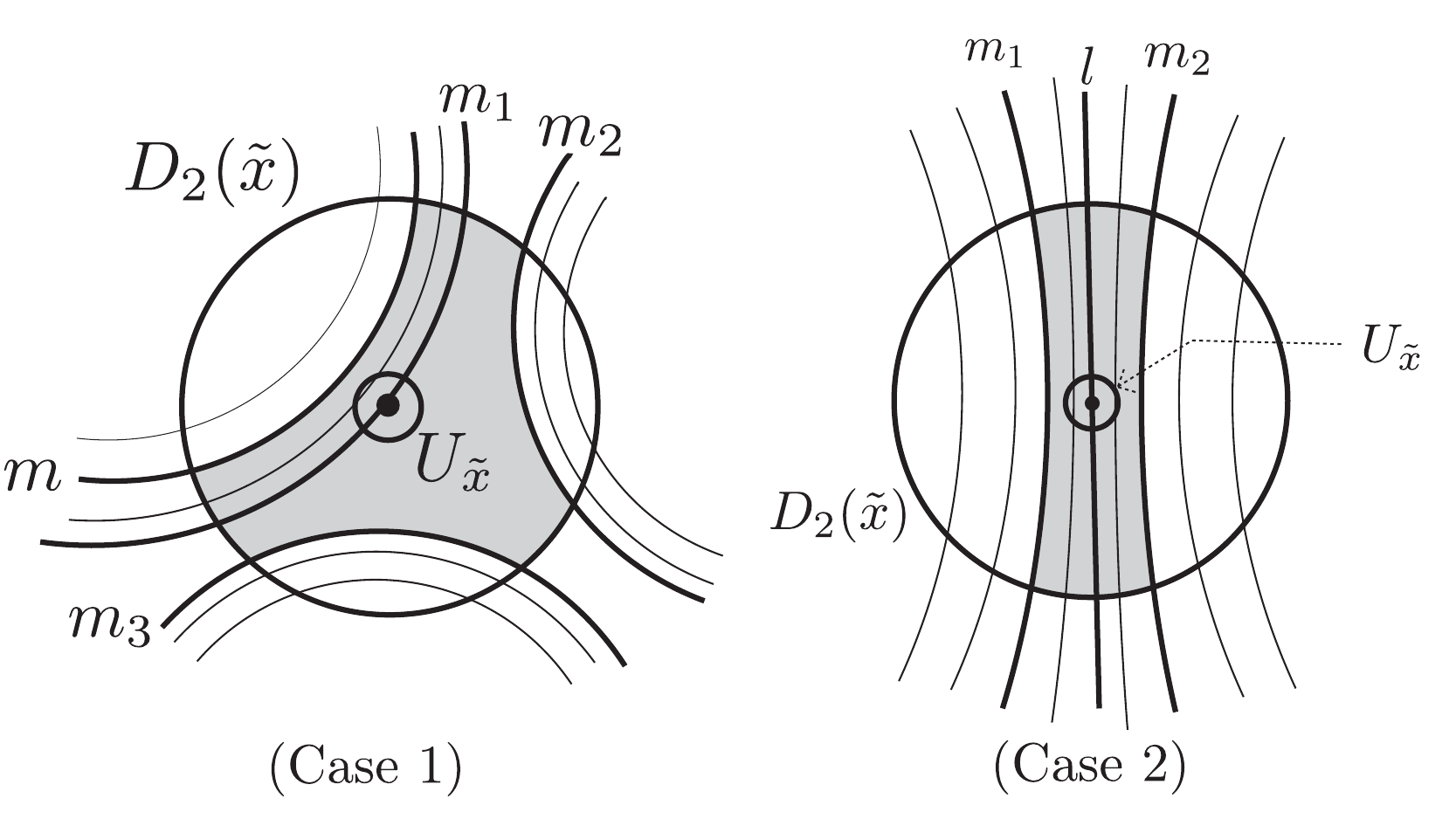}
\caption{In each picture, $P \cap D_2(\td{x})$ is shaded.}
\label{Case1&2}
\end{figure}

\begin{proof}
We claim that, for every $x \in (S, \tau)$, there exist a neighborhood $U_x$ of $x$  and $i_x \in \N$ such that, if $i > i_x$, then $\omega(s) < \epsilon$ for every $s \in \mathcal{NT}(l_i)$  with $s \cap U_x \neq \emptyset$.   
This would imply the Lemma, since $S$ is compact.
Let $\td{M} = (\td{\nu}, \td{\omega})$ and $\td{l}_i$ denote the total lifts of $M$ and $l_i$ to $\h^2$, respectively. 
Choose a lift $\td{x}$ of $x$ to $\h^2$. 
Through the covering map from $\h^2$ to $(S, \tau)$, the above claim is equivalent to the following:
There exist a neighborhood $U_{\td{x}}$ of $\td{x}$ and $i_{\td{x}} \in \N$ such that, if $i > i_{\td{x}}$, then $\td{\omega}(s) < \epsilon$ for every $s \in \mathcal{NT}(\td{l}_i)$ such that $s \cap U_{\td{x}} \neq \emptyset$.

{\it Case 0.} Suppose first that $x \not\in |\nu|$; Then, $\td{x} \not\in |\td{\nu}|$. 
Let $P$ be the component of $\h^2 \setminus |\td{\nu}|$ that contains $\td{x}$. 
Then $P$ is an open convex region bounded by some leaves of $\td{\nu}$. 
Clearly, only finitely many such boundary leaves intersect $D_{2}(\td{x})$. 
Let $m_1, m_2, \ldots , m_k$ denote these intersecting leaves. 
Let $m'_h = D_2(\td{x}) \cap m_h$ for all $h \in \{1,2, \ldots, k\}$. (See Figure \ref{PandV}.) 
Recall that there is a projection $\mc{P}: \h^2 \to T$,  where $T$ is the tree dual to $\td{M}$ (\S \ref{DT}).
By Lemma \ref{box-nhd}, for every $\epsilon > 0$ and every $h \in \{1, 2, \ldots, k\}$, there exists a flow box neighborhood 
$V_h$ of  $m'_h$ that projects to a geodesic segment of  length less than $\epsilon/2$ in $T$.
Let $s_h$ denote the geodesic segment $\mc{P}(V_h)$.
Since $m_h$ is a boundary geodesic of $P$, the leaves of $\td{M}$ contained in $P$ do not accumulate to $m_h$. 
Therefore, by the construction of $V_h$, we can assume that  $\mc{P}(m_h') = \mc{P}(m_h) = \mc{P}(P)$ is an end point of $s_h$.                                                                                                                                                                                                                                                                                                                                                                                                                                                                                                                                                                                                                                                                                                                                                                                                                     
Therefore, $\cup_h V_h$ projects to $\vee_h s_h \st T$, the one point union of $s_h$ that identifies the end points $\mc{P}(m_h)$ of $s_h$. 
Then $~(\cup_h V_h) \cup (P \cap D_2(\td{x}))$ also projects onto $\vee_h s_h$,
	and the diameter of $\vee_h s_h$ is less than $\epsilon$. 

Since the sequence $(\td{l_i})$ approximates $\td{\nu}$, each $m_h$ is approximated by a sequence $(n_{h,i})_i$ such that $n_{h,i}$  is a leaf of  $\td{l}_i$. 
For the rest of Case 0, we always assume that $i$ is sufficiently large. 
For each $i$, let $P_i$ be the open convex region bounded by  $\bigsqcup_h n_{h,i}\>$.
Then $(P_i)$ limits to $P$ as $i$ goes to infinity. 
Hence, $~\td{x} \in P_i \cap D_2(\td{x}) \subset  (\cup_h V_h) \cup ( P \cap D_2(\td{x}))$. 
Take an open neighborhood $U_{\td{x}}$ of $\td{x}$ such that the closure of $U_{\td{x}}$ is contained in the interior of  $P \cap D_1(\td{x})$. 
Then $U_{\td{x}}$ is contained in $P_i$. 
Since $\pt P_i \st |\td{l}_i|$,  
every $s \in \mc{NT}(\td{l}_i)$ with $s \cap U_{\td{x}} \neq \emptyset$ is contained in $cl(P_i) \cap \mr{D}_2(\td{x})$, and therefore $s$ is contained in  $(\cup_h V_h) \cup ( P \cap D_2(\td{x}))$.
Hence $\mc{P}(s) \st \vee_h s_h$ and $\td{\omega}(s) < \epsilon$. 
\begin{figure}[htbp]
\includegraphics[width = 5in]{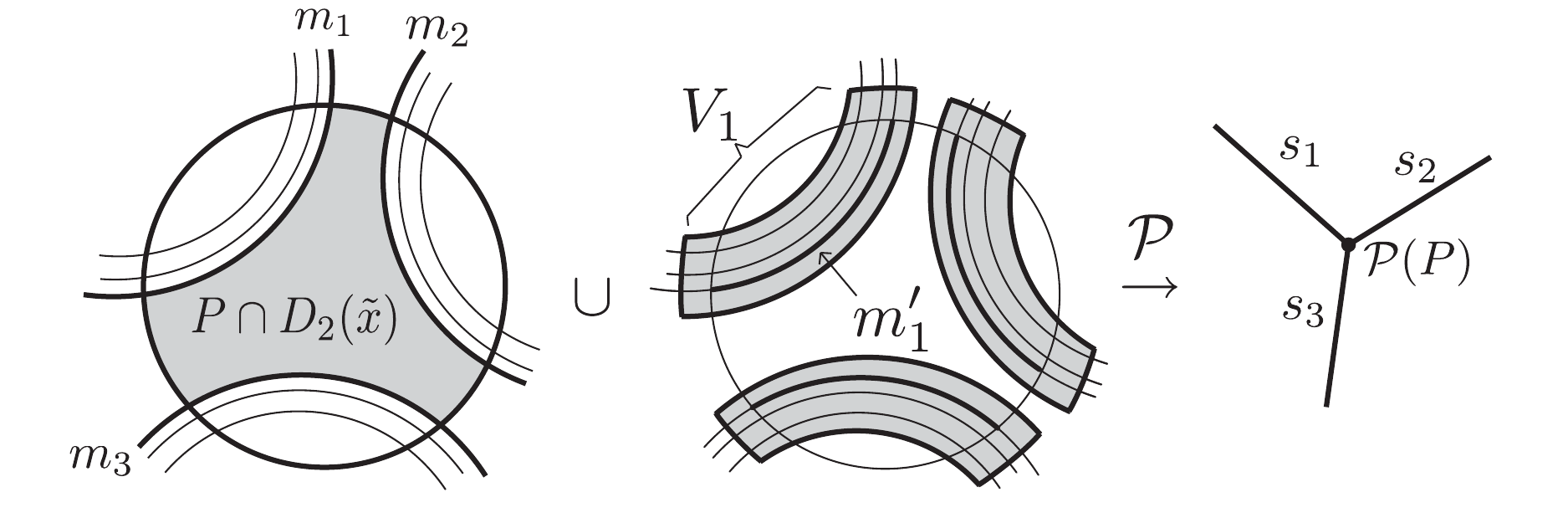}
\caption{Case 0.}
\label{PandV}
\end{figure}

Next, suppose $x \in |\nu|$; then, $\td{x} \in |\td{\nu}|$.
Let $l \in \td{\nu}$ be  the leaf containing $\td{x}$, and let ${l'} = l \cap D_2(\td{x})$.
There are two ways that $\td{\nu}$ can accumulate to $l$:\\
{\it Case 1.} the leaves of $\td{\nu}$ accumulate from only one side of $l$, or\\
{\it Case 2.} the leaves of $\td{\nu}$ accumulate from both sides of $l$. \\
In both cases, by Lemma \ref{box-nhd}, there exists a flow box neighborhood $V$ of $l'$ with height less than $\epsilon/2$. 

In Case 1, $~l$ is a boundary geodesic of a component $P'$ of $\h^2 \setminus |\td{\nu}|$. 
Let $m_1, m_2, \ldots, m_k$ be the boundary geodesics of $P'$ that intersect $D_2(\td{x})$. 
We can assume that $m_1 = l$. 
Then $l$ bounds a half plane disjoint from  $P'$, and the leaves of $\td{\nu}$ contained in this half plane accumulate to $l$.
Consider a  leaf $m$ of $\td{\nu}$ in this half plane, such that $m \cap D_2(x) \st V$. 
Let $P$ be the open convex region in $\h^2$ bounded by $m$ and $m_2, m_3, \ldots, m_h$, which contains $\td{x}$. 
Then $\mc{P}(P)$ is a geodesic segment of length less than $\epsilon/2$.
Take  an open neighborhood $U_{\td{x}}$ of $\td{x}$ such that the closure of $U_{\td{x}}$ is contained in the interior of $P \cap D_1(\td{x})$.
The same argument as Case 0 shows that $U_{\td{x}}$ satisfies the desired property. (See the left picture in Figure \ref{Case1&2}.) 

In Case 2, in each component of $\h^2 \sm l$, consider a leaf $m_i$, $i =1, 2$, of $\td{\nu}$ close to $l$.
Let $P$ be the open convex region in $\h^2$ bounded by $m_1$ and $m_2$, which contains $l$.
By taking $m_1$ and $m_2$ sufficiently close to $l$, we can assume that $P \cap D_2(\td{x})$ is contained in $V$. 
Take an open neighborhood $U_{\td{x}}$ of $\td{x}$ such that the closure of $U_{\td{x}}$ is contained in the interior of  $P \cap D_1(\td{x})$. 
Again, as in Case 0, we see that $U_{\td{x}}$ satisfies the desired property. (See the right picture in Figure \ref{Case1&2}.)  
\end{proof}

Let $Q$ be a component of $S \sm l_i$, and let $\td{Q}$ be a lift of $Q$ to $\h^2$. Note that, if we take a different lift $\td{Q}$, then $I(\td{M}, \td{Q})$ changes only by an element of $\pi_1(S)$.  
In particular, $\|I( \td{M}, \td{Q})\|$ does not depend on the choice of $\td{Q}$.

 \begin{proposition}\label{PropId} 
For every $\epsilon > 0$, there exists $i_0 \in \N$ such that, if $i > i_0$,
then $\|I( \td{M}, \td{Q})\| < \epsilon$ for every component $Q$ of $S \sm l_i$.
\end{proposition}

\begin{figure}[htbp]
\includegraphics[width= 2in]{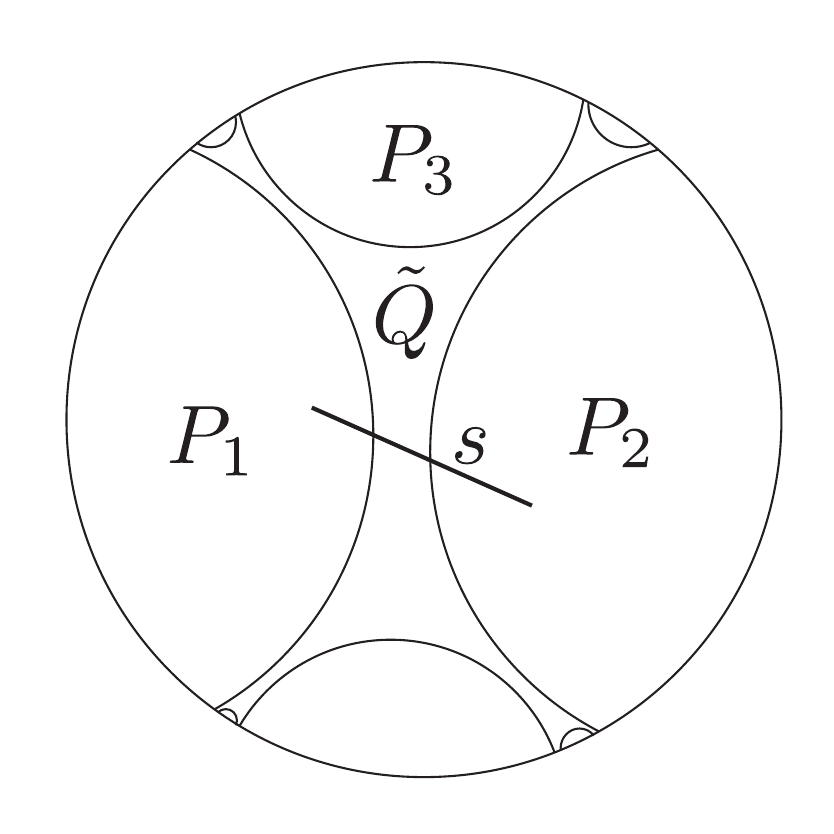}
\caption{}
\label{prop7-4}
\end{figure}

\begin{proof}
For every $Q$, its lift $\td{Q}$ is an open convex region bounded by some leaves of $\td{l}_i$.
Therefore, $\h^2 \sm \td{Q}$ consists of closed half planes bounded by these leaves. 
We can assume that $\td{Q}$ is $\pi_1(Q)$-invariant.  
Let $k$ be the number of the boundary components of $Q$.
The $\pi_1(Q)$-action permutes the complementary half planes, and this action on the half planes has exactly $k$ orbits.
Let $H_1, H_2, \ldots, H_k$ be the representatives of the orbits.
Since $l$ consists of $n$ disjoint simple loops on $S$, $~k \leq 2n$. 

Let $H$ be a component of $\h^2 \sm \td{Q}$. 
Note that $\pt H$ is transversal to $\td{M}$.
Then a leaf $l$ of $I(\td{M}, \td{Q})$ intersects $H$ if and only if $l$ intersects $\pt H$. 
Therefore, $ I(I(\td{M}, \td{Q}), H) = I(\td{M}, \pt H)$.
By the covering map from $\h^2$ to $(S, \tau)$, $~\pt H$ covers some $l_{i,j}$.
Recalling that $\tl_{i,j}$ is a lift of $l_{i,j}$ to $\h^2$, $~I(\td{M}, \pt H) \cong I(\td{M}, \tl_{i,j})$.
Since $I(\td{M}, \td{Q})$ and  $I(I(\td{M}, \td{Q}), H)$ coincide on $H$,
 $$\| I(\td{M}, \td{Q})|_H \| = \| I(I(\td{M}, \td{Q}), H)|_H \| \leq \| I(I(\td{M}, \td{Q}), H) \| = \| I(\td{M}, \tl_{i,j}) \|.$$
Since $\td{Q}$ is convex, for every geodesic segment $s$ in $\h^2$,  $~\td{Q} \cap s$ is either empty or a geodesic segment.
In addition, $s \cap (\h^2 \sm Q)$ consists of at most 2 geodesic segments, each of which is contained in a component of $\h^2 \sm Q$ (see Figure \ref{prop7-4}).  
Therefore, by the definition of the norm,
 \begin{eqnarray}
 &&    \|I( \td{M}, \td{Q})\| \leq   \|I( \td{M}, \td{Q})|_{\td{Q}} \| + 2\max\{\| I( \td{M}, \td{Q})|_{H_j}  \|~ | ~j = 1, \ldots, m\} \nonumber \\
 &&     \hspace{0.8in} \leq   \|I( \td{M}, \td{Q})|_{\td{Q}} \| + 2\max\{\| I( \td{M}, \td{l}_{i,j})  \|~ | ~j = 1, \ldots, m\}. \label{eq}
\end{eqnarray}

By Lemma \ref{LemId},  for every $\epsilon > 0$, there exists $i_0 \in \N$ such that, if $i > i_0$, then $\|I( \td{M}, \td{Q})|_{\td{Q}} \| < \epsilon$ for every component $Q$ of $S \sm l_i$. 
By Proposition \ref{propA},  for every $j$, $~ \| I( \td{M}, \td{l}_{i,j}) \| \to 0$ as $i \to \infty$.
Therefore, for every $\epsilon > 0$, if $i$ is sufficiently large, then (\ref{eq}) is bounded from above by $\epsilon$ for every $Q$. 
\end{proof}

\begin{proof}(Theorem \ref{id})
Recall that $m$ is the number of the periodic leaves of $L$. 
For each $k \in \{1, 2, \ldots, m\}$, let $A_k = \kappa^{-1}(p_k) \st (S,C)$, which is a flat cylinder of height $w(p_k)$, and let $A = \bigsqcup _{k= 1}^m A_k$.
Since $S \sm A$ is homeomorphic to $S \sm p$ by $\kappa$,  $~\kappa^{-1}(l_{i,j})$ is a simple closed loop on $(S, C)$ for each sufficiently large $i \in \N$ and for each $j \in \{1, 2, \ldots, m\}$.
Let $l'_{i,j} = \kappa^{-1}(l_{i,j})$ and $l_i' = \bigsqcup_j {l'}_{i,j}$. 
Then $\kappa$ restricts to a homeomorphism from $S \setminus (A \bigsqcup l_i')$ to $S \sm (p \bigsqcup l_i)$. 
In particular, the components of  $S \setminus (A \bigsqcup l_i')$ and  $S \sm (p \bigsqcup l_i)$ correspond bijectively. 
Note that, since $p_h$'s and $l_{i,j}$'s are not pairwise homotopic (when $i$ is fixed), each component of  $S \setminus (A \bigsqcup l_i')$ and  $S \sm (p \bigsqcup l_i)$ is (the interior of) a compact surface with boundary of negative Euler characteristic. 

Let $S'_i = S \setminus (A \bigsqcup l_i')$. 
First, we shall show that, when $i$ is large, every component $R'$ of $S'_i$ is admissible.
Let $R = \kappa(R')$, which is a component of $(S, \tau) \sm (l_i \bigsqcup p)$.
Let  $\td{R}$ and $\td{R'}$ be the corresponding lifts of $R$ and $R'$ to the universal covers,  $(\td{S}, \td{C})$ and $(\h^2, \td{L})$, respectively, so that $\td{\kappa}(\td{R'}) = \td{R}$.
Choose $\delta > 0$ that satisfies Theorem \ref{ib}. 
Let $I_R = I(\td{M}, \td{R} )$.
Then, since $R \cap p = \emptyset$, $~I_R = I(\td{L}, \td{R})$.
Note that $R$ is contained in a component of $S \sm l_i$.
Therefore, by Proposition \ref{PropId}, $~\| I_R \| < \delta$ for every component $R'$ of $S'_i~$, provided that $i$ is sufficiently large. 
Hence, by Corollary \ref{admissibleD}, the projective structure $C(I_R)$ on $\oD$ is admissible.
Since $I_R = I(\td{L}, \td{R})$ and $\td{R}$ is a convex subset of $\h^2$ bounded by geodesics, by Lemma  \ref{C(I)/stC(L)}, $~C(I_R)$ embeds into $C(\td{L})$.
Each boundary leaf of $R$ does not transversally intersect a leaf of $L$ with positive weight and $R$ is not a geodesic. 
Therefore, by Corollary \ref{BasicCoro}, $~\mc{R}(\td{C}, \td{R}') = \mc{R}( C(I_R), \kappa_{I_R}^{-1}(\td{R})) \st C(I_R)$, where $\kappa_{I_R}$ is the collapsing map for $C(I_R)$.
Since $C(I_R)$ is admissible,  $\mc{R}(\td{C}, \td{R}')$ is also admissible.

Through its action, $\pi_1(R) \cong \pi_1(R')$ is regarded as a Schottky group in $\psl$.
Let $\beta_{I_R}: \h^2 \to \h^3$ be the bending map induced by $I_R$.
Then, by Theorem \ref{ib}, $~\beta_{I_R}$ is an injective quasiisometric embedding, and it extends to an equivariant embedding of $\pt \h^2$  to $\RS$.
In particular, this extension takes the limit set of $\pi_1(R)$ to the limit set of $\rho(\pi_1(R))$ homeomorphically and $~\rho_{\pi_1(R)}$-equivariantly. 
Therefore, $\rho|_{\pi_1(R)}$ is an isomorphism onto a Schottky group in $\PSL$.
Hence, the restriction of $C$ to $R'$ is admissible. 


We have given a desired decomposition of $(S, C)$, except that the flat cylinders $A_k$ are not integral. 
In what follows, instead of cutting out the whole $A_k$  from $S$, we cut out a maximal integral flat cylinder contained in $A_k$.
Taking the union of the maximal integral flat cylinders and $l'_i$, we shall show that the complementary regions of the union are admissible, which completes the decomposition.

For  $x \in \R_{\geq 0}$, let $[x] = \max \{ n \in \mathbb{Z}_{\geq 0} \>|\> 2\pi n \leq x\}$.
Besides, let $$a_k = \frac{w(p_k) - 2\pi[w(p_k)]}{ 2}$$ for each $k$. 
Then $0 \leq a_k < \pi$. 
(\emph{Remark}: 
One might wonder why we only have $a_k < \pi$, instead of having $a_k < \pi/2$ as in Assumption (iv) of Theorem \ref{ib2}.
The following example illustrates a ``hidden'' $\pi/2$-annulus, which fills in this difference. 
Assume that $\Gamma$ is a Schottky group in $\psl$.
Let $H$ be the convex hull, in $\h^2$, of the limit set of $\Gamma$. 
Then $H/\Gamma$ is a projective structure on a  surface with boundary. 
Then the projective structure $\h^2/\Gamma$ can be obtained by the attachment of a $\pi/2$-annulus to each boundary component of $H/\Gamma$.)

For each $k \in \{1,2, \ldots, m\}$, we cut each $A_k$ along two admissible loops into three flat cylinders of heights  $a_k$, $2\pi[w(p_k)]$, $a_k$ in this order.
The middle cylinder $A'_k$ is integral and the others are not.
(See $A_1$ in Figure \ref{PandR}.)
If $w(a_k) < 2\pi$, then $A'_k$ degenerates to an admissible loop in the middle of $A_k$ (see $A_2$ in Figure \ref{PandR}).  

Let $A' = \bigsqcup_{k = 1}^m A'_k$ and $S''_i = S \setminus (l_i' \bigsqcup A')$. 
Each non-integral flat cylinder obtained above shares exactly one boundary component with a component of $S'_i = S \sm (l'_i \bigsqcup A)$.
Therefore, each component $R'$ of $S'_i$ is contained in a component $P'$ of $S''_i$. 
If a boundary circle $l$ of $R'$ maps to a periodic leaf $p_k$  via $\kappa$, then $l$ bounds a flat cylinder of height $a_k$ in $P'$. 
If $l$ maps to an approximating loop $l_{i,j}$ by $\kappa$, then $l$ is a boundary component of $P'$.
Thus,  each component $P'$ of $S_i''$ is the union of a component  $R'$ of $S'_i$ and the non-integral flat cylinders sharing a boundary component with $R'$. 
In particular, $R'$ is a deformation retract of $P'$. 
Letting  $B_1, B_2, \ldots, B_r$ be these non-integral flat cylinders, set $$~P' = R' ~\cup~ ( B_1 \cup \ldots \cup B_r)$$ (the shaded region in Figure \ref{PandR}).
Then $p_{b_1}:=\kappa(B_1)$, $p_{b_2}:=\kappa(B_2),$ $~\ldots,~ p_{b_r}:=\kappa(B_r)$ are  periodic leaves of $L$.

We have $\kappa(R') = \kp(P') =: R$.
Let $\td{R}' \st \td{P}'$ be lifts of $R'$ and $P'$ to $\td{S}$, respectively.
Then $\td{\kp}(\td{R}') = \td{\kp}(\td{P}') =: \td{R}$, which is a lift of $R$ to $\h^2$.
Let $\lambda_\partial$ be the geodesic lamination on $\h^2$ consisting of the boundary geodesics of $\td{R}$ that are lifts of periodic leaves of $L$. 
Let $\theta = max\{ a_k - \pi/2,~ 0 ~| ~k = 1, 2, \ldots, m \}$; then $0 \leq \theta < \pi/2$. 
Assign the weight $\theta$ to each leaf of $\lambda_\pt$, and obtain  a measured lamination $L_\pt$ on $\h^2$.
Recall that $I_R = I(\td{L}, \td{R}) = I(\td{M}, \td{R})$.
Since there are no leaves of $L$ intersecting both $\td{R}$ and $|L_\pt|$, and $L_\pt$ consists of isolated leaves of $\td{L}$,
	therefore, $~|L_\pt|$ and $|I_R|$ are disjoint.
Then let $L_P = I_R \bigsqcup L_\partial$.
Each leaf $l$ of $L_\pt$ is a boundary geodesic of $\td{R}$, and  each leaf of $I_R$ intersects $\td{R}$ but does not intersect leaves of $L_\pt$.
Therefore, each $l$ is an outermost leaf of $L_P$.

\begin{figure}[htbp]
\includegraphics[width = 5in]{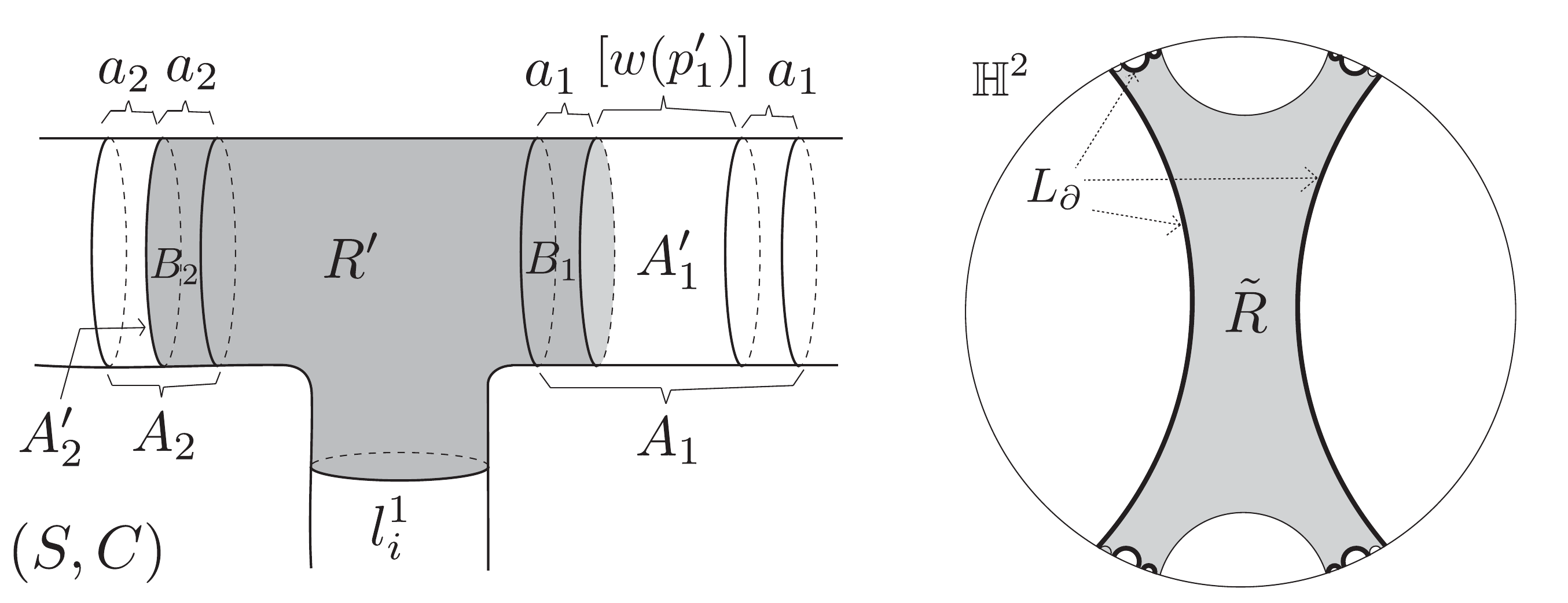}
\caption{In the left picture, the region $P'$ is shaded. In the right picture, the bold lines are the leaves of $L_P$.}
\label{PandR}
\end{figure}

We now apply Theorem \ref{ib2} with $L = L_P$ and $L' = L_\pt$.
We just have checked that $L_\pt \st \pt L$ (Hypothesis (i) of Theorem \ref{ib2}).
Every leaf of $L_\pt$ has the weight $\theta < \pi/2$ (Hypothesis (iv)).
Since $p_1, \ldots, p_m$ are disjoint simple loops on $S$ and $L_\pt$ consists of lifts of these loops, there exists $D > 0$ such that every pair of distinct leaves of $L_\partial$ has a  distance greater than $D$ (Hypothesis (iii)). 
Apply Theorem \ref{ib2} to this $D$, and obtain $\delta, S, T > 0$.
Apply Proposition \ref{PropId} with $\epsilon = \delta$, and obtain $i_0 \in \N$. 
Let $Q$ be the component of $S \sm l_i$ containing $R$.
Then $I(\td{L}, \td{R}) = I(\td{M}, \td{R}) \st I (\td{M}, \td{Q})$, where $\td{Q}$ is a lift of $Q$ to $\h^2$ (that contains $\td{R}$). 
Therefore, by the Proposition \ref{PropId}, if $i > i_0$, then $\| I_R \| \leq \| I(\td{M}, \td{Q})\| < \delta$ (Hypothesis (ii)). 
By Theorem \ref{ib2}, if $i$ is sufficiently large, $~{\beta}_{L_P}$ is an injective $(S, T)$-quasiisometric embedding for every component $P$ of $S''_i$.

By Corollary \ref{admissibleD}, the projective structure $C(L_P)$ on $\oD$ is admissible.
Let $\kappa_{I_R}$ and $\kappa_{L_P}$ be the collapsing maps associated with $C(I_R)$ and $C(L_P)$, respectively. 
Since $I_R = I(L_P, \td{R}) = I(L, \td{R})$, by Corollary \ref{BasicCoro},  we have $\mc{R}(\td{C},  \td{R}') = \mc{R}(C(I_R),  \kappa_{I_R}^{-1}(\td{R})) = \mc{R}( C(L_P), \kappa_{L_P}^{-1}(\td{R}) )$.
Therefore, $\mc{R}(\td{C},  \td{R}')$ is admissible. 

For $h \in \{1, 2,\ldots, r\}$, let $\td{p}_h$ be a lift of $p_{b_h}$ to $\h^2$ that bounds $\td{R}$.
Let $H$ be the component of $\h^2 \sm \td{R}$ bounded by $\td{p}_h$.
Then, since $\td{R}$ is open, $H$ is a closed half plane.
Note that $H$ is uniquely determined by the choice of $h$ and the choice of the lift of $p_{b_h}$.  
Observe that $|L_P| \cap H = \td{p}_h$. 
Then $\mc{R}( C(L_P), \kappa_{L_P}^{-1}(H \sm \td{p}_h)) $ is a crescent of angle $\pi/2$,  and $\mc{R}( C(L_P), \kappa_{L_P}^{-1}(\td{p}_h))$ is a crescent of angle $\theta$.
Therefore, $\mc{R}( C(L_P), \kappa_{L_P}^{-1}(H) )$ is a crescent of angle $\pi/2 + \theta$, and it is a component of  $C(L_P) \sm \mc{R}( C(L_P),  \kappa_{L_P}^{-1}(\td{R}) )$.

Each component of $\mc{R}(\td{C}, \td{P}')\sm \mc{R}(\td{C}, \td{R}')$ is a lift $\td{B}_h$ of some $B_h$ to $\td{S}$.
There is a lift $\td{p}'_h$ of $p_h$ separating $\td{B}_h$ and $\td{R}'$ in $(\td{S}, \td{C})$ (Figure \ref{RandB}).
The height of $\td{B}_h$ is $a_h \leq \pi/2 + \theta$.
By the argument above, when $\mc{R}(\td{C}, \td{R}') = \mc{R}( C(L_P), \kappa_{L_P}^{-1}(\td{R}) )$ is embedded in $C(L_P)$, $~\td{p}'_h$ bounds a component of $C(L_P) \sm \mc{R}(\td{C}, \td{R}')$, which is a crescent of angle $\pi/2 + \theta$. 
Therefore, this embedding extends to the embedding of $\td{B}_h \cup \mc{R}(\td{C}, \td{R}')$ to $C(L_P)$ (see Figure \ref{RandB}).
Different components of $\mc{R}(\td{C}, \td{P}') \sm \mc{R}(\td{C}, \td{R})$ correspond to different $\td{B}_h$.
Therefore,  the embedding extends disjointly  to all components of $\mc{R}(\td{C}, \td{P}') \sm \mc{R}(\td{C}, \td{R})$ and we obtain an embedding of $\mc{R}(\td{C}, \td{P}')$ into $C(L_P)$.
By the construction, the embedding of $\mc{R}(\td{C}, \td{R})$ is $\pi_1(R')$-equivariant. 
Since $\mc{R}(\td{C}, \td{P}')$ embeds in $C(L_P)$, and $C(L_P)$ has an injective developing map, therefore $\mc{R}(C, P')$ has an injective developing map. 

\begin{figure}[htbp]
\includegraphics[width = 3in]{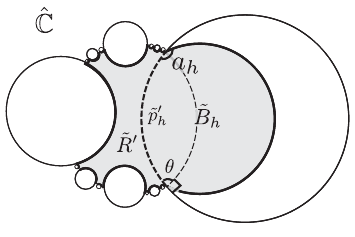}
\caption{$\mc{R}(\td{C}, \td{R}') = \td{R}'$ and $\td{B}_h$ in $C(L_P)$}
\label{RandB}
\end{figure}

Since $R'$ is a deformation retract of $P'$, $~\pi_1(R')$ is equal to $\pi_1(P')$ as subgroups of $\pi_1(S)$.
In particular, $\rho|_{\pi_1(P')} = \rho|_{\pi_1(R')}$.
We have already seen that  $\pi_1(R')$ is isomorphic to  a purely loxodromic subgroup of $\PSL$ via $\rho$.
Therefore, so is $\rho_{\pi_1(P')}\>$. 
Hence, the restriction of $C$ to $P'$ is admissible.  
\end{proof}


\end{document}